\documentclass[12pt]{article}
\scrollmode

\usepackage{color}

\usepackage{amsfonts}
\usepackage{amsmath}
\usepackage{amsthm}
\usepackage{amssymb}
\usepackage{graphicx}

\usepackage{color}

\allowdisplaybreaks

\newtheorem{Theo}{Theorem}%[section]

\newtheorem{lem}[Theo]{Lemma}%[theo]
\newtheorem{prop}[Theo]{Proposition}%[theo]
\newtheorem{rmk}[Theo]{Remark}%[section]
\newtheorem{cor}[Theo]{Corollary}%[theo]%[section]

\def\Z{{\mathbb Z}}
\def\R{{\mathbb R}}
\def\P{{\mathbb P}}
\def\E{{\mathbb E}}
\def\H{{\mathbb H}}

\def\cp{{\cal P}}

\def\cs{{\cal S}}
\def\cq{{\cal Q}}
\def\ch{{\cal H}}

\def\ci{{\cal I}}
\def\cw{{\cal W}}
\def\car{{\cal R}}
\def\ct{{\cal T}}
\def\cj{{\cal J}}
\def\cn{{\cal N}}
\def\ce{{\cal E}}

\def\d{\delta}

\def\={&=&}
\def\eps{\epsilon}

\def\lf{\lfloor}
\def\rc{\rceil}
\def\rf{\rfloor}
\def\lc{\lceil}
\def\mr{\mathring}
\def\le{\left}
\def\ri{\right}

\def\s{\sigma}
\def\a{\alpha}
\def\b{\beta}
\def\t{\theta}

\def\o{\omega}
\def\g{\gamma}
\def\G{\Gamma}
\def\l{\lambda}
\def\La{\Lambda}
\def\vf{\varphi}
\def\k{\kappa}

\def\nn{\nonumber}

\def\ds{\displaystyle}

\renewcommand{\cal}{\mathcal}

\begin{document}

\title{Scaling limit of the radial Poissonian web}

\date{}

\author{L.~R.~G.~Fontes %\footnote{Corresponding author} 
\footnote{IME-USP, Rua do Mat\~ao 1010, 05508-090,
S\~ao Paulo SP,  Brazil, lrfontes@usp.br} \thanks{Partially
supported by CNPq grant 305760/2010-6, and FAPESP grant 2009/52379-8}
\and
L.~A.~Valencia H.\footnote{Depto.~de Matem\'aticas, UdeA, 
Calle 67 No.~53-108, Medell\'\i n, Colombia, lavalench@gmail.com} 
%\thanks{????}
\and
G.~Valle\footnote{Instituto de Matemática, Universidade Federal do Rio de Janeiro, Caixa Postal 68530, 21945-970, 
Rio de Janeiro RJ, Brasil, glauco.valle@im.ufrj.br}
\thanks{Partially supported by CNPq grants 304593/2012-5 and 474233/2012-0}
}
\maketitle

\begin{abstract}
We consider a variant of the radial spanning tree introduced by Baccelli and Bordenave. 
Like the original model, our model is a tree rooted at the origin, built on the realization of a planar Poisson point process. 
Unlike it, the paths of our model have independent jumps.
We show that locally our diffusively rescaled tree, seen as the collection of the paths connecting its sites to the root, 
converges in distribution to the {\em Brownian Bridge Web}, which is roughly speaking a collection of coalescing Brownian bridges 
starting from all the points of a planar strip perpendicular to the time axis, and ending at the origin.  
\end{abstract}

%%%%%%%%%%%%%%%%%%%%%%%%%%%%%%%%%%%%%%%%%%%%%%%%%%%%%%%%%%%%%%%%%%%%%%%%%%%%%%%%%%
%%%%%%%%%%%%%%%%%%%%%%%%%%%%% INTRO %%%%%%%%%%%%%%%%%%%%%%%%%%%%%%%%%%%%%%%%%%%%%%
%%%%%%%%%%%%%%%%%%%%%%%%%%%%%%%%%%%%%%%%%%%%%%%%%%%%%%%%%%%%%%%%%%%%%%%%%%%%%%%%%%

%
\section{Introduction}
\label{intro}

\setcounter{equation}{0}

The Radial Spanning Tree (RST) of a Poisson point process has been introduced by Baccelli and Bordenave in \cite{bb}. 
It is a random planar tree rooted at the origin whose vertices are points in a realization of a homogeneous Poisson point 
process on the plane. The motivation comes from an increasing interest in random graphs from both applied and 
theoretical fields. For instance, spanning trees are an essential modeling tool in communication networks ---
see \cite{bb} and references in this respect. 
From a theoretical point of view, the RST is related to models like the random minimal directed 
spanning trees~\cite{br,pw1,pw2} and the Poisson trees~\cite{flt,ffw}. 

Consider the branches of the RST as random paths heading towards the root of the tree. 
It forms a system of planar colescing random paths heading towards the origin. A large class of space-time 
systems of one dimensional symmetric coalescing random walks converge in distribution under diffusive scaling 
to the Brownian Web (BW), which is a system of coalescing Brownian motions introduced in~\cite{a}, and later 
studied in~\cite{tw} and~\cite{finr1}. The latter paper started the study of convergence of rescaled systems 
of random paths to the BW (see also~\cite{finr}), and since then quite a few papers have presented convergence results to the BW and 
variations~\cite{bmsv,cfd,cv,cv1,ffw,finr2,fn,nrs,ss,ss1,rss}, to name a few, 
of which~\cite{ffw} is concerned with the convergence of (the rescaled paths of) a Poisson tree to the BW.

Two other of these papers are worth singling out, since they are are closely related either to the present paper 
or to~\cite{bb}. It is worth to first point out two features of the paths of the radial spanning tree of~\cite{bb}, 
distinguishing them from the case of most of the previously studied models. One feature is the (long-range) 
dependence of their increments on past history, leading to non-Markovianness, and the other is their radial 
character: they are directed towards the origin. 

\cite{rss} deals with a planar model of discrete space-time (non-radial) paths with long-range dependence
of increments similar to those of the radial spanning tree, and shows that when diffusively scaled, they converge
to the BW. A key step in it showing that the dependence washes out in the limit is by establishing the 
existence of regeneration times with well-behaved tails when the path increments loose memory (a similar step is
taken in the analysis in~\cite{bb}, even though that model is better described as a continuum space-time one,
and the regeneration times there have a somewhat weaker character).

\cite{cv1} deals with a radial path model in contiunuum space, discrete time, whose paths have independent 
increments, showing convergence of the rescaled paths to a web of coalescing Brownian bridges, and object
called Radial Brownian Web in that paper --- let us call it the Brownian Bridge Web (BBW) in this paper. 
It was obtained as a suitable mapping of the ordinary BW. 
(The appearance of coalescing Brownian bridges is somehow to be expected in such a situation, and 
even in the context of the RST, even though it has not been suggested before, as far as we know.)

%Considering systems of coalescing random paths, for Poisson trees convergence to the Brownian Web have been 
%shown in~\cite{ffw}. 

Based on these results, it is natural to ask whether the (suitably rescaled) radial spanning tree converges 
(locally) in distribution to the BBW. 
%Related to this question, Coletti and Valencia \cite{cv1} have also 
%introduced a radial spanning tree rooted at the origin, but vertices in their tree are on circles of 
%integer radius centered at the origin. They call it the discrete radial Poisson tree and prove its convergence 
%in distribution under diffusive scaling to a continuous transformation of the Brownian Web named 
%the Radial Brownian Web. 

The aim of this paper is to introduce, as a variation of the RST, a path model which we call the 
{\em Radial Poissonian Web}. Like the RST, it consists of radial
paths directed towards the origin, passing through the points of a planar Poisson point process, and whose increments, 
unlike the case of the RST, have independent (step) increments (from one Poisson point to the next).
We show that such a model, when diffusively rescaled around a ray, converges to the BBW. 

Our approach is similar to that of~\cite{cv1}. We roughly speaking transform the model (from the beginning, in our case) 
to a planar, non-radial path model, verify converge criteria of the latter model to the BW, and map back, thus giving
rise to the BBW. 

We believe that our model is considerably closer to the RST than the one of~\cite{cv1}, and poses extra technical 
issues in its analysis. Given the results of~\cite{rss}, and similar if weaker ones of~\cite{bb} (concerning 
regeneration times), it seems safe to {\em conjecture} that an approach similar to ours here as far as transforming
the radial model to a planar one is concerned, and similar to the one of~\cite{rss} in the analysis of the planar model,
by resorting to regeneration times (like the ones already to shown exist in~\cite{bb}, probably with some modification),
will lead to a proof that the BBW is the scaling limit of the RST, thus answering a question asked in~\cite{bb}.

Some other features of our model (to be described in detail in the next section) and our approach, with respect to
others in the literature are as follows.
The random paths in our system have (long-range) dependence one to another before coalescence. 
%and exhibit long-range dependence 
%(as far as fixed deterministic times are concerned, since one has to look back from a deterministic 
%time along a path for the last Poisson path it passed through to know the law of the next step; the step by step increments,
%as described above, are independent). 
%These last two properties are 
This feature is already present in drainage network models whose convergence 
to the Brownian Web were considered in \cite{cfd,cv}. The techniques employed in the study of convergence in these 
two last cited papers and also in \cite{nrs} are going to play a central role in our approach, although~\cite{cv} 
and~\cite{nrs} deals with systems with crossing paths, which is not the present case, as we will see in detail below.

%So far we are not able to prove that the scaling limit of Baccelli and Bordenave's tree is the Radial Brownian Web or a 
%variation of it. The problem is related to how the radial spanning tree is built. Applying to the radial spanning tree 
%the transformation that maps the radial Brownian Web onto the standard Brownian Web, we obtain a system of space-time 
%coalescing random paths which have non Markovian jumps. 
%A given point $x$ in the Poisson point process will be connect by an edge to the closest point between those that are closer than $x$ to the origin. 

The paper is organized as follows: In Section 2 we introduce our model, the RPW, and main result, namely its convergence 
in distribution when suitably centered and scaled to the BBW. In Section~\ref{proof}, we delineate the strategy of proof of 
our result, describing the mapping from the radial to the planar model. Sections~\ref{pflm} and~\ref{pfpr} are devoted to 
preliminaries, leading to a listing of the criteria of convergence of the transformed, planar model to the BW.
In Section~\ref{sec:coaltime} we prove an auxiliary result about the coalescence time of two of our paths, and in the
remainder Sections~\ref{sec:I} and~\ref{sec:B1E} we verify the aforementioned convergence criteria.

%%%%%%%%%%%%%%%%%%%%%%%%%%%%%%%%%%%%%%%%%%%%%%%%%%%%%%%%%%%%%%%%%%%%%%%%%%%%%%%%%%
%%%%%%%%%%%%%%%%%%%%%%%%%%%%% MOD %%%%%%%%%%%%%%%%%%%%%%%%%%%%%%%%%%%%%%%%%%%%%%
%%%%%%%%%%%%%%%%%%%%%%%%%%%%%%%%%%%%%%%%%%%%%%%%%%%%%%%%%%%%%%%%%%%%%%%%%%%%%%%%%%

\section{The model. Main result.}\label{mod}

\setcounter{equation}{0}

Let $\cp$ be a planar Poisson point process of intensity 1 and let $\cp_0=\cp\cup\{O\}$, 
where $O$ denotes the origin. 
We form a network of paths in the plane {\em directed} towards the origin as follows. 
From each point $x$ in $\cp$ we start a directed edge ending at another 
point $s=s(x)$, the {\em successor} of $x$, in $\cp_0$, chosen as follows. 

\begin{figure}
 \begin{center}
\input{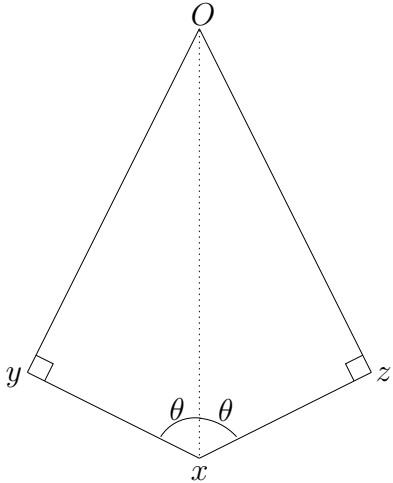}
\end{center}
 \caption{Representation of $\cq_x$ (in full straight lines).}
\label{fig:qx}
\end{figure}

\paragraph{Choice of successor}
First consider the convex quadrangle 
$\cq_x:=\cq_{x,\theta}$ with opposing vertices at the origin and at $x$. 
Let $y=y(x)$ and $z=z(x)$ be the other two vertices of $\cq_x$. We assume that the internal angles at
$y$ and $z$ are right ones. The internal angle at $x$ is given by $2\t$, where $\t\in(0,\pi/2)$ is a
parameter of the model. We further assume that the segment $\overline{Ox}$ bisects the angles at the 
origin and at $x$. See Figure~\ref{fig:qx}.

We  now  choose $s$ as the point of $\cp_0$ within $\cq_x$ which is farthest from the origin. 
This choice is uniquely defined for every $x\in\cp$ for a.e.~$\cp$.
(A perhaps simpler and more natural choice of $s$ would be as the point of $\cp_0\cap\cs_x$ closest to $x$, 
where $\cs_x$ is the circular sector centered at $x$ whose circumference 
contains the origin, with central angle $2\t$ bisected by the segment $\overline{Ox}$;
but it poses difficulties which we have rather avoid in this paper. 
We believe our results can be established also for this choice at 
least for $\t$ small enough --- $\t<\pi/4$ should do ---, with not much more work in the latter case.) 

Since $\cq_x$ lies within the circumference centered at the origin passing by $x$, 
we have that $s$ is (a.s.) closer to the origin than $x$ is.

\paragraph{Paths}
For each $x\in\cp$, let us now introduce the sequence $\g_x=(s_n)=(s_n(x))_{n=0}^I$ as follows: 
$s_0(x)=x$ and $s_i=s(s_{i-1}(x))$, $i=1,\ldots,I$, where $I=I(x)$ is such that $s_i\ne O$, $i=0,\ldots I-1$,
and $s_I=O$. It is a straightforward matter to check that $I$ is almost surely well defined and finite.
Then $\g_x$ may be identified with a directed path starting at $x$ and ending at the origin,
passing through the edges $e_i=(s_{i-1},s_i)$, $i=1,\ldots,I$. We will have $\g_x$ actually as the planar 
(polygonal) curve determined by $\{s_i,\,i=0,\ldots,I\}$ by linear interpolating between the edge endpoints, 
in the usual way. We say that $x$ is the starting point of $\g_x$, with the origin its ending point.

%We note that from our assumption on the angle $\t$, $\g_x,\,x\in\cp$, are almost surely well defined.

Then $\G=\{\g_x,\,x\in\cp\}$ is a family of paths from every point of $\cp$ ending at the origin.
We want to understand the large scale behavior of $\G$ under {\em diffusive} scaling. We will then
indeed consider a (relatively small) portion of $\G$, whose paths will indeed be also clipped 
at a point when they get (macroscopically) close to the origin, and differently modified
if they wander too far to the sides, in a way to be explained in detail below.

\paragraph{Modified paths}
In order to define which portion of $\G$ and which clipping and other modification we will consider, 
let us denote a general point $x$ in the plane by its polar coordinates, namely, in complex numbers 
notation, $x=re^{i\vf}$.

Let $\a\in(0,1)$ and $1/4<a<b<1/2$. Let us first define the clipping. For $x$ such that
$r\in[\a n,n]$ and $\vf\in[-\pi/2\pm n^{-b}]$ (with the notation $[c\pm d]:=[c-d,c+d]$ for real numbers $c,d$).
Let $I'=\min\{0\leq i\leq I:\|s_i(x)\|<\a n\}$, where $\|\cdot\|$ is the Euclidean norm in the plane.
Finally, let $\g_x'$ be the path determined by $\{s_i(x),\,i=0,\ldots,I'-1\}$ as a planar polygonal curve, 
similarly as in the definition of $\g_x$, concatenated to the segment 
$(s_{I'-1}(x),s'_{I'}(x))$, where $s'_{I'}(x)$ is the intermediate point of the segment
$(s_{I'-1}(x),s_{I'}(x))$ such that $\|s'_{I'}(x)\|=\a n$.
%where that segment intersects the circumference centered at the origin with radius $\a n$.
%$(s_{I'-1}(x),\a e^{i\arg(s_{I'-1}(x))})$.
So $\g_x'$ is $\g_x$ clipped at the point where it is at distance $\a n$ from the origin. 

Now the final modification: again for $x$ such that $r\in[\a n,n]$ and $\vf\in[-\pi/2\pm n^{-b}]$, let 
$I''=\min\{0\leq i\leq I':|\arg(s_i(x))+\pi/2|>n^{-a}\}\wedge I'$, where $\min\emptyset=\infty$; 
then $\g_x''$ is the path determined by $\{s_i(x),\,i=0,\ldots,I''-1\}$ concatenated to the segment
$(s_{I''-1}(x),\a e^{i\arg(s_{I''-1}(x))})$.

\begin{figure}
\begin{center}
\input{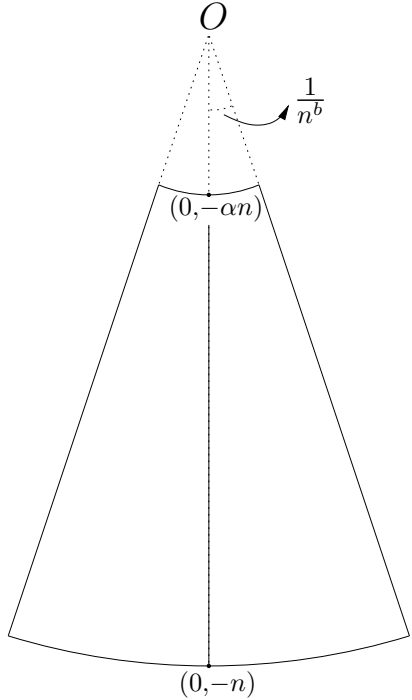}
\end{center}
\caption{Representation of $\La_n$ (in full lines).}
\label{fig:ln}
\end{figure}

The set of paths we will analyse is then as follows.
\begin{equation}
 \label{paths}
 \G_n=\{\g''_x,\,,x\in\cp\cap\La_n\},
\end{equation} 
where
\begin{equation}
 \label{lan}
 \La_n=\La_n(\a,b)=\{x=re^{i\vf},\,r\in[\a n,n],\,\vf\in[-\pi/2\pm n^{-b}]\}.
\end{equation} 
See Figure~\ref{fig:ln}.

Note that as subsets of the plane, the paths of $\G_n$ are contained in 
\begin{equation}
 \label{lanp}
 \bar\La_n:=\La_n(\a,a).
\end{equation} 

We will show below that {\em with high probability} (that is, with probability going to 1 as $n\to\infty$)
we have $\g''_x=\g'_x$ for all paths $\g''_x$ in $\G_n$. $\G_n$ has %almost sure 
properties, to be discussed below, which are convenient for our analysis.

\paragraph{Diffusive rescaling of $\G_n$}

\begin{figure}
 \begin{center}
\input{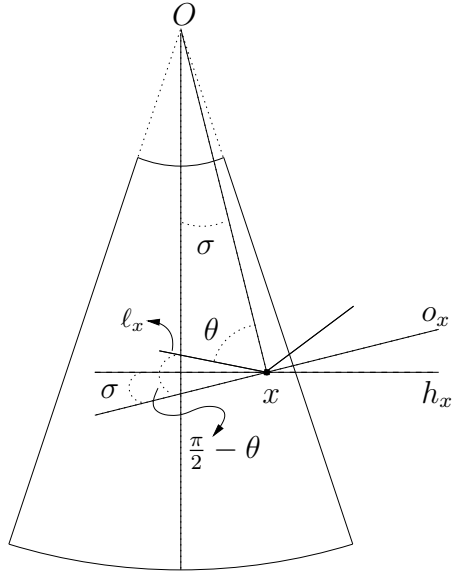}
\end{center}
 \caption{Illustration of the fact that $\cq_x$ (whose bottom sides appear partially in full lines;
 see Figure~\ref{fig:qx}) is above the horizontal (dashed) line through $x\in\bar\La_n$ for all $n$ large enough.}
\label{fig:pathprop}
\end{figure}

We first claim that for all $n$ large enough and every $x\in\bar\La_n$, $\cq_x$ is above the horizontal
line through $x$. To argue this, let us consider $x\in\bar\La_n$ such that $\arg(x)\geq-\pi/2$.  
By the definition of $\bar\La_n$, we have that $\s:=\arg(x)+\pi/2\leq n^{-a}$. Let $h_x$ and $o_x$ be the 
horizontal and perpendicular straight lines through $x$, respectively. We thus have that the angle between 
$h_x$ and $o_x$ to the left of $x$ equals $\s$. Since the angle between $\overline{Ox}$ 
and the left bottom side of $\cq_x$, call it $\ell_x$, equals $\t$, and by assumption $\t<\pi/2$,
we have that the angle between $\ell_x$ and $h_x$ to the left of $x$ equals $\frac\pi2-\t-\s$, which is
positive for all large enough $n$, and thus $\ell_x$ and $\cq_x$ are above $h_x$, and the claim 
is established. See Figure~\ref{fig:pathprop}.

We may then conclude that %for a.e.~$\cp$, 
all paths of $\G_n$ have their second coordinate strictly increasing
as one goes from their starting points to the origin, and we may then interpret the second coordinate as
time and the paths as continuous trajectories.

Let us now introduce a {\em diffusive} scaling of the plane given by the map
$R:\R^2\to\R^2$ such that  
\begin{equation}
 \label{eq:tt}
 D(x_1,x_2)=\left(\frac{x_1}{\sqrt n},\frac{x_2}n\right).
\end{equation} 
%(thinking of the second coordinate as time, as pointed out above). 
The rescaled set of trajectories is then
\begin{equation}
 \label{rescpaths}
 \tilde\G_n=\{D(\g);\,\g\in\G_n\},
\end{equation} 
where $\tilde\g=D(\g)$ is the image of $\g\in\G_n$ by $D$, and may be readily checked 
to be a trajectory. % in $\Pi_{-1}^{-\a}$.

We want to state a (weak) convergence result for $\tilde\G_n$ in terms of a (limiting) random set of
trajectories in a suitable path space. Let us introduce this space now.

\paragraph{Path space; Hausdorff space}

Let $\b,\b',\b''$ be real numbers such that $\b<\b'\leq\b''$, and for $\b\leq t_0\leq\b'\leq t_1\leq\b''$, 
let $C[t_0,t_1]$ denote the set of functions $f$ from
$[t_0,t_1]$ to $[-\infty,\infty]$ such that $\tanh\circ f$ is continuous, and define
\begin{eqnarray}
 \label{pathspace}
 \Pi=\Pi_{\b}^{\b',\b''}=\bigcup_{\b\leq t_0\leq\b'\leq t_1\leq\b''}C[t_0,t_1]\times\{(t_0,t_1)\},
\end{eqnarray} 
where $(f,t_0,t_1)\in\Pi$ represents a path/trajectory in $[-\infty,\infty]\times[\b,\b'']$ starting at
$(f(t_0),t_0)$ and ending at $(f(t_1),t_1)$. Note that, given $\a'\in(0,\a)$, we have that 
$\tilde\G_n\subset\Pi_{-1}^{-\a,-\a'}$ for all large enough $n$.

For $(f,t_0,t_1)\in\Pi$, let $f^\star:[\b,\b'']\to[-\infty,\infty]$ be the function which
equals $f$ in $[t_0,t_1]$, identically equals $f(t_0)$ in $[\b,t_0]$, and identically equals $f(t_1)$ in $[t_1,\b'']$.
Let us equip $\Pi$ with the metric
\begin{eqnarray}
 \label{metric}
 d((f,t_0,t_1),(g,s_0,s_1))=|t_0-s_0|\vee|t_1-s_1|\vee\sup_{\b\leq s\leq\b''}|\tanh(f^\star(s))-\tanh(g^\star(s))|,
\end{eqnarray} 
under which it is complete and separable. One may check that $\tilde\G_n$ is a Borel subset of $\Pi_{-1}^{-\a,-\a'}$ 
(as soon as $n$ is large enough; $\a'\in(0,\a)$, as above).

Consider now the Hausdorff space $\ch=\ch_{\b}^{\b',\b''}$ of compact subsets of $(\Pi,d)$ equipped with the Hausdorff metric
\begin{eqnarray}
 \label{hmetric}
 d_{\ch}(K,K')=\sup_{h\in K}\inf_{h'\in K'}d(h,h')\vee\sup_{h'\in K'}\inf_{h\in K}d(h,h'),
\end{eqnarray} 
which makes it a complete separable metric space.

The weak limit of $\tilde\G_n$ will be given below in terms of the {\em restricted} Brownian web, an object
which we describe next.

\paragraph{Restricted Brownian web}

The (ordinary) Brownian web is (the closure of) a family of coalescing one dimensional Brownian motions 
(ordinarily with diffusion coefficient 1, but here we will for convenience take the diffusion coefficient
a number $\o$, to be exhibited below --- see~(\ref{eq:om}) ---, which is a function of $\t$, 
with $\o\in(0,\infty)$ whenever $\t\in(0,\pi/2)$) starting from all points in the planar space-time. 
See [FINR] for details.
In this paper we will consider a {\em restricted version} of that family: we will take Brownian paths 
(with diffusion coefficient $\o$) starting 
from times {\em in the interval} $[0,\tau]$ and {\em ending at time} $\tau$, where $\tau=\frac1\a-1$. 
One way to define/analyse this object is to construct/study it as done in [FINR], but taking 
$\Pi_{0}^{\tau}:=\Pi_{0}^{\tau,\tau}$ and $\ch_{0}^{\tau}:=\ch_{0}^{\tau,\tau}$ as the relevant path and 
sample spaces, instead of the full space versions considered in [FINR], which we will call here
$(\bar\Pi,\bar d)$ and $(\bar\ch,\bar d_{\bar\ch})$ (in [FINR], the full space versions were denoted 
as $(\Pi,d)$ and $(\ch,d_{\ch})$).

A more economic approach though is to take a {\em restriction map} of the full space (ordinary) Brownian web.
We discuss that now.

Let $\bar\cw$ be the ordinary Brownian web, defined in the Hausdorff space $\bar\ch$ of compact 
subsets of $\bar\Pi$ (see [FINR] for details; again, beware of the different notation: $\bar\Pi$ and $\bar\ch$
here correspond to $\Pi$ and $\ch$ in [FINR], respectively). 

Let $\bar\Pi^\tau$ be the (closed) subset of paths
of $\bar\Pi$ starting at or before time $\tau$, and let $R:\bar\Pi^\tau\to\Pi_{0}^{\tau}$ be the restriction
of semi-infinite paths of $\bar\Pi^\tau$ to $[0,\tau]$, namely, given $(f,t_0)\in\bar\Pi^\tau$, in which case
$f:[t_0,\infty]\to[-\infty,\infty]$ is continuous (under a suitable metric on $[-\infty,\infty]^2$ making it
compact --- see [FINR] for details), and $t_0\leq\tau$, 
\begin{equation}
\label{rest1}
R((f,t_0))=(f',t_0,\tau),%\in\Pi_{0}^{\tau},
\end{equation} 
where $f'=f|_{[t_0,\tau]}$ is the restriction of $f$ to $[t_0,\tau]$. 

Note that the metric induced by $\bar d$ on $\Pi_{0}^{\tau}$ via $R$ is equivalent to $d$.
It is a straightforward matter to check that the map $R$ is continuous.

Similarly, let $\bar\ch^\tau$ be the (closed) subset of compact sets of $\bar\ch$ whose paths belong to 
$\bar\Pi^\tau$ (in other words, whose paths start at or before time $\tau$),
and let $\car:\bar\ch^\tau\to\ch_0^\tau$ be such that 
\begin{equation}
\label{rest2}
\car(K)=\{R(h);\,h\in K\}.
\end{equation} 

Note that the metric induced by $\bar d_{\bar\ch}$ on $\ch_{0}^{\tau}$ via $\car$ is equivalent to 
$d_{\ch_{0}^{\tau}}$.
It is a straightforward matter to check that the map $\car$ is continuous.

Let $\bar\cw^\tau=\bar\cw|_{\bar\Pi^\tau}$ be the restriction of the ordinary Brownian web to sets 
of paths of $\bar\Pi^\tau$ (that is, to sets of paths starting at or before $\tau$).
We define the {\em restricted} Brownian web as 
\begin{equation}
\label{rbw}
\cw^\tau_0:=\car(\bar\cw^\tau).
\end{equation} 
In other words, $\cw^\tau_0$ is the restriction of the ordinary (full space) Brownian web to 
subsets of paths starting between times $0$ and $\tau$, and ending/clipped at $\tau$.

\paragraph{Main result}

In order to state our main result concerning $\tilde\G_n$, let us introduce the following maps. 
Let $\psi: [-\infty,\infty]\times[0,\tau]\to[-\infty,\infty]\times[-1,-\a]$
such that 
\begin{equation}
 \label{psi}
 \psi(y,s)=\left(\frac y{1+s},-\frac1{1+s}\right).
\end{equation}
Next let $\psi': \Pi_0^\tau\to\Pi_{-1}^{-\a}$
be such that given a path $(f,s_0,s_1)\in\Pi_0^\tau$, 
\begin{equation}
 \label{psi1}
 \psi'((f,s_0,s_1))=(\psi\circ f,\psi(s_0),\psi(s_1)).
\end{equation}
$\psi'((f,s_0,s_1))$ may be also described as the image by $\psi$ of the path $(f,s_0,s_1)$
as a set in $[-\infty,\infty]\times[0,\tau]$.
Finally, let $\psi'':\ch_0^\tau\to\ch_{-1}^{-\a}$
be such that given $K\in\ch_0^\tau$, 
\begin{equation}
 \label{psi2}
 \psi''(K)=\{\psi'(h),\,h\in K\}.
\end{equation}

We are now ready to state our main result.

\begin{Theo}\label{main}
 Let $\a'\in(0,\a)$ be given. As $n\to\infty$,
 \begin{equation}
  \label{eq:main}
  \tilde\G_n\Rightarrow\psi''(\cw^\tau_0),
 \end{equation}
 where $\cw^\tau_0$ is the restricted Brownian web given in~(\ref{rbw}) above, and ``$\Rightarrow$'' 
 stands for convergence in distribution in $\ch_{-1}^{-\a,-\a'}$.
\end{Theo}

Notice that $\psi''(\cw^\tau_0)\in\ch_{-1}^{-\a,-\a}$, which is a closed subset of $\subset\ch_{-1}^{-\a,-\a'}$.

It is a straightforward matter to verify that the path of $\psi''(\cw^\tau_0)$ starting at a deterministic 
point $(y,s)\in(-\infty,\infty)\times[0,\tau]$ (from the well known corresponding property of the (restricted) 
Brownian web, there is almost surely only one such path) is a Brownian bridge (with diffusion coefficient $\o$)
starting at $(y,s)$ and finishing
at the origin, stopped at time $\tau$. For this reason we may call $\psi''(\cw^\tau_0)$ the {\em Brownian bridge
web}, which then may be roughly described as a collection of coalescing Brownian bridges starting from all points
of $[-\infty,\infty]\times[0,\tau]$ and finishing at the origin, stopped at time $\tau$.

%%%%%%%%%%%%%%%%%%%%%%%%%%%%%%%%%%%%%%%%%%%%%%%%%%%%%%%%%%%%%%%%%%%%%%%%%%%%%%%%%%
%%%%%%%%%%%%%%%%%%%%%%%%%%%%% PROOF %%%%%%%%%%%%%%%%%%%%%%%%%%%%%%%%%%%%%%%%%%%%%%
%%%%%%%%%%%%%%%%%%%%%%%%%%%%%%%%%%%%%%%%%%%%%%%%%%%%%%%%%%%%%%%%%%%%%%%%%%%%%%%%%%

\section{Proof of Theorem~\ref{main}}\label{proof}

\setcounter{equation}{0}

%\subsection{Preliminaries}\label{pre}

\paragraph{Strategy.}

In this section we (begin to) present a proof of Theorem~\ref{main}. The rough idea is to map 
$\tilde\G_n$ with the inverse of $\psi$ to a set of paths of $\Pi_0^{\tau}$, prove a convergence 
result of the mapped set to the restricted Brownian web $\cw^\tau_0$, and then map back with $\psi$.

%Due to a technical issue
For convenience, we will actually apply a variation of this strategy. We will take the 
rescaling of $\G_n$ after a suitable mapping --- related to, but not quite the inverse of $\psi$ ---, 
in this case to $\Pi_0^{\tau n}$, and then prove two things:
1) the convergence to the restricted Brownian web, and 2) that the map by $\psi$ of the rescaled image
of $\G_n$ is close to $\tilde\G_n$. 
%The second point happens to be straightforward; the main body of work is thus point 1. 

\paragraph{Map.}

Representing points of $\bar\La_n$ in complex polar coordinates, namely $x = re^{i(−\frac\pi2+\s)}$ , where
$r = |x|$ and $\s = \frac\pi2 + \arg(x)$, with $\a n \leq r \leq n$ and $|\s| \leq n^{-a}$ , let 
$\Xi  : \bar\La_n \to [−\infty, \infty] \times [0, \tau n]$
such that
\begin{equation}\label{eq:map}
 \Xi  (re^{i(−\frac\pi2 +\s)}) = \left(n\s,\frac{n-r}{r/n}\right),
\end{equation} 
and let
\begin{equation}\label{eq:gp}
 \G'_n=\{\Xi (\g):\,\g\in\G_n\},
\end{equation} 
where, given a path $\g\in\G_n$, $\g'=\Xi  (\g)$ is its image by $\Xi $, which happens to be a path in $\R\times[0,\tau n]$.

Notice that 
\begin{eqnarray}
 \label{lp}
 &\La'_n:=\Xi (\La_n)=[-n^{1-b},n^{1-b}]\times[0,\tau n];&\\
 \label{blp}
 &\bar\La'_n:=\Xi (\bar\La_n)=[-n^{1-a},n^{1-a}]\times[0,\tau n].&
 \end{eqnarray}

Let us now rescale $\G'_n$ diffusively. Let
\begin{equation}
 \label{rescppaths}
 \tilde\G'_n=\{D (\g');\,\g'\in\G'_n\},
\end{equation} 
where $D$ was given in~(\ref{eq:tt}) above, and
$\tilde\g'=D(\g')$ is the image of $\g'\in\G'_n$ by $D$, and may be readily checked 
to be a trajectory in $\Pi_{0}^{\tau}$.

The proof of Theorem~\ref{main} then follows readily from the following two auxiliary results.

\begin{prop}\label{aux}
 As $n\to\infty$,
 \begin{equation}
  \label{eq:aux}
  \tilde\G'_n\Rightarrow\cw^\tau_0,
 \end{equation}
 where $\cw^\tau_0$ is the restricted Brownian web given in~(\ref{rbw}) above, and ``$\Rightarrow$'' 
 stands for convergence in distribution.
\end{prop}

\begin{lem}\label{lm:dist}
Let $\a'\in(0,\a)$ be given. We have
\begin{equation}
\label{eq:dist}
d_{\ch_{-1}^{-\a,\a'}}(\tilde\G_n,\psi''(\tilde\G'_n))\to0
\end{equation} 
with high probability as $n\to\infty$.
\end{lem}

In the next section, we prove Lemma~\ref{lm:dist} and begin the proof of Proposition~\ref{aux}, deferring 
the conclusion to the remaining sections of the paper.
%The proofs of Proposition~\ref{aux} and Lemma~\ref{lm:dist} will be deferred to the next section.
%In that section, we will indeed, again for convenience, modify $\G'_n$ a further time before 
%proving Proposition~\ref{aux}.

%, involving in the first case verifying standard criteria for convergence to the Brownian web.

%%%%%%%%%%%%%%%%%%%%%%%%%%%%%%%%%%%%%%%%%%%%%%%%%%%%%%%%%%%%%%%%%%%%%%%%%%%%%%%%%%
%%%%%%%%%%%%%%%%%%%%%%%%%%%%% PROOF %%%%%%%%%%%%%%%%%%%%%%%%%%%%%%%%%%%%%%%%%%%%%%
%%%%%%%%%%%%%%%%%%%%%%%%%%%%%%%%%%%%%%%%%%%%%%%%%%%%%%%%%%%%%%%%%%%%%%%%%%%%%%%%%%

\section{Proofs of Proposition~\ref{aux} and Lemma~\ref{lm:dist}}
\label{proofaux}

\setcounter{equation}{0}

\subsection{Preliminaries}\label{pre}

We will indeed, also for convenience, work with a variant of $\G_n$, which equals $\G_n$ with high 
probability. Let us go into this point next, and after that describe properties of the map $\Xi $
and their consequences for our analysis.

\paragraph{Variant of $\G_n$.}

For $x\in\La_n'$ and $0\leq l\leq\log n$,  let $w=w(x,l)$ be the point on the segment $\overline{Ox}$ 
at distance $l$ from $x$, and let $\ct_{x,l}$ denote the closed subset of $\cq_x$ outside the circumference 
centered at the origin and passing by $w$. See Figure~\ref{fig:qcirc}. 
%below the straight line orthogonal to 
%$\overline{Ox}$ which passes by $w$. %One readily checks that for all 
%outside the open circle centered at the origin whose circumference passes by $w$. 
\begin{figure}
 \begin{center}
\input{qcirc.pstex_t}
\end{center}
 \caption{Representation of $\ct_{x,l}$ (in full lines). The length of the dotted segment inside $\ct_{x,l}$ 
 is $l$, as indicated.}
\label{fig:qcirc}
\end{figure}
Note that for all $n$ large enough, and all $x\in\La_n'$, $\ct_{x,l}$ is roughly 
triangular, as shown in the picture, and not roughly pentagonal, which could happen if we took $w$ at distance 
of order $n$ from $x$ (to visualize the latter point, consider a high enough point $w$ of $\overline{Ox}$ in
Figure~\ref{fig:qx}.

We define now a set of polygonal paths $\hat\G_n=\{\hat\g_x,\,x\in\cp\cap\La_n\}$ starting in $\La_n$ made of 
the concatenation edges similarly as before, namely for $x\in\cp\cap\La_n$, $\hat\g_x$ is determined by 
$\{\hat s_i(x),\,i=0,\ldots,\hat I\}$ , where $\hat s_0(x)=x$, and given $\hat s_{i-1}(x)\in\bar\La_n$, 
if $\ct_{\hat s_{i-1}(x),\log n}\cap\cp\ne\emptyset$, then $\hat s_{i}(x)$ is the point of 
$\ct_{\hat s_{i-1}(x),\log n}$ which is farthest from the origin; if $\ct_{\hat s_{i-1}(x),\log n}\cap\cp=\emptyset$, 
then $\hat s_{i}(x)$ is given by $w(\hat s_{i-1}(x),\log n)$;
%, with $w(\cdot,\cdot)$ defined at the beginning of the latter paragraph
and $\hat I=\min\{i>0:\hat s_{i}(x)\notin\bar\La_n\}$.
Then $\hat\g_x$ is the path within $\bar\La_n$ obtained by linear interpolating $\hat s_i(x)$, $i=0,\ldots,\hat I-1$, 
concatenated to an edge from $\hat s_{\hat I-1}(x)$, say $\hat e$, defined as follows. Let $\hat s'_{\hat I}(x)$ be 
the point of the edge $(\hat s_{\hat I-1}(x),\hat s_{\hat I}(x))$ where that edge intersects the boundary of $\bar\La_n$. 
Then $\hat s'_{\hat I}(x)$ is either at the top of that bounday, or it is at its sides. If it is at the top, then 
$\hat e=(\hat s_{\hat I-1}(x),\hat s'_{\hat I}(x))$; otherwise, 
$\hat e=(\hat s_{\hat I-1}(x),\a e^{i\arg(\hat s_{\hat I-1}(x))})$.

We will now show that with high probability, $\hat\G_n=\G_n$.

\begin{lem}\label{lm:id}
With high probability
\begin{equation}
\label{eq:id}
\hat\g_x=\g''_x\,\mbox{ for all }x\in\cp\cap\bar\La_n.
\end{equation} 
\end{lem}

\noindent{\bf Proof} 

It is enough to show that with high probability $s_1(x)=\hat s_1(x)$ for all $x\in\cp\cap\bar\La_n$,
but that follows immediately from proving that with high probability  
$\ct_{x,\log n}\cap\cp\ne\emptyset$ for all $x\in\cp\cap\bar\La_n$.

This in turn follows readily from the following two estimates. First, that outside an event of vanishing
probability as $n\to\infty$ the cardinality of $\cp\cap\bar\La_n$ is bounded above by constant
times $n^{2-a}$ (since that is the order of magnitude of the area of $\bar\La_n$. And the second estimate is
for the probability of $\ct_{x,\log n}\cap\cp=\emptyset$ for a single $x\in\ct_{x,\log n}\cap\cp=\emptyset$, 
which equals $e^{-\mbox{area of $\ct_{x,\log n}$}}\leq e^{-\mbox{area of $\ct'_{x}$}}$, where $\ct'_x$ is the triangle 
obtained by removing the portion of $\ct_{x,\log n}$ above the line through $w(x,\log n)$ orthogonal to $\overline{Ox}$. 
Since the area of $\ct'_x$  is given by $\sin\t(\log n)^2$ , the second
result follows from the fact that $n^{2-a}e^{-\sin\t(\log n)^2}\to0$ as $n\to\infty$. $\square$

\medskip

We may then work with either $\G_n$ or $\hat\G_n$ . We will for a while below work with $\hat\G_n$. We will
next discuss some properties of the map $\Xi $.

From the rules of formation of paths of $\hat\G_n$, we may describe those of paths of 
\begin{equation}\label{eq:hgp}
\hat\G'_n=\{\Xi (\hat\g);\,\hat\g\in\hat\G_n\}, 
\end{equation}
where $\hat\g'=\Xi (\hat\g)$ is the image of $\hat\g\in\hat\G_n$ 
(which as an immediate corollary to Lemma~\ref{lm:id} equals $\G'_n$ with high probability).
First let us look at $\cp'$, the image of $\cp$ by $\Xi $.

\begin{lem}\label{lm:pp}
$\cp'$ is a Poisson point process on $\R\times[0, \tau n]$ with intensity measure
\begin{equation}
\label{eq:int}
\frac1{(1+s/n)^3}\,dy\,ds,
\end{equation} 
where $s$ is the second (time) coordinate.
\end{lem}

\noindent{\bf Proof} 

Follows from the distribution of $\cp$ and a straightforward computation of the Jacobian of
the relevant transformation. $\square$

\medskip

\begin{figure}
 \begin{center}
\input{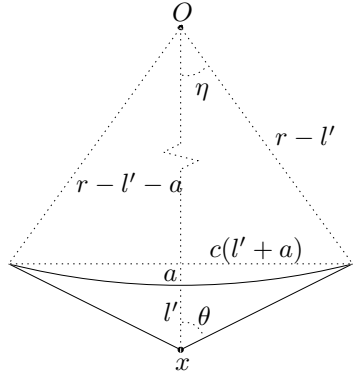}
\end{center}
 \caption{$\eta = \arcsin\frac{c(l' +a)}{r-l'}$; $a$ satisfies $c^2(l' + a)^2 + (r - l' - a)^2 = (r - l' )^2$ }
\label{fig:tt}
\end{figure}

\begin{figure}
 \begin{center}
\input{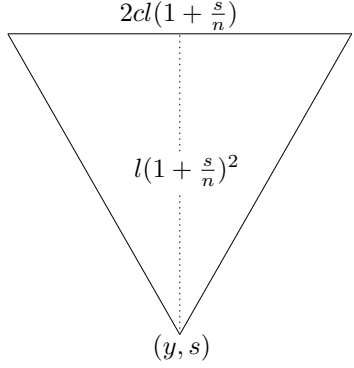}
\end{center}
 \caption{Rough shape of $\ct'_{(y,s),l}$. Top side of triangle is orthogonal to $s$-axis.}
\label{fig:rshape}
\end{figure}

We will next see how the choice of successors of the points determining the paths of $\hat\G_n$ get
translated by $\Xi $. Given Lemma~\ref{lm:pp}, a key step in that direction is understanding 
how $\ct_{x,l}$, $x\in\bar\La_n$, $0 \leq l \leq \log n$, are
mapped by $\Xi $. A straightforward analysis finds that for $(y, s)\in\R^2$ such that $\Xi  (x) = (y, s)$
\begin{equation}\label{eq:xt}
 \Xi (\ct_{x,l})=\ct'_{(y,s),l}=
 \bigcup_{l'\in[0,l]}[y\pm n\eta]\times
 \left\{s+\frac{(1+\frac sn)^2}
             {1-\frac{l'}n\left(1+\frac sn\right)}l'\right\},
\end{equation} 
where 
\begin{equation}\label{eq:eta}
\eta = \arcsin\frac{c(l' +a)}{r-l'} 
=\arcsin\left\{\frac cn\frac{1+\frac sn}{1-\frac{l'}n\left(1+\frac sn\right)}(l'+a)\right\},
\end{equation} 
with $c=\tan\t$, $r = n/(1 + s/n)$, and $a$ satisfies 
\begin{equation}\label{eq:a}
c^2(l' + a)^2 + (r - l' - a)^2 = (r - l' )^2. 
\end{equation} 
See Figure~\ref{fig:tt}: notice that the top vertices of $\ct_{x,l'}$ in that picture get mapped
to $\{y - n\eta, y + n\eta\} \times\left\{s+\frac{{(1+\frac sn)^2}}{1-\frac{l'}n\left(1+\frac sn\right)}l'\right\}$.
Notice that for $n$ large, $\ct'_{(y,s),l}$ is roughly the isosceles
triangle depicted on Figure~\ref{fig:rshape}.

\paragraph{Description of $\hat\G'_n$.} We may now describe the paths of $\hat\G'_n$ in terms of starting points and
successors in $\cp'$ similarly as we did above with $\hat\G_n$, using $\ct'_{(y,s),l}$ instead of $\ct_{x,l}$. 
It is an important step to understand the succession mechanism for the image of $\{\hat s_i(x);\,i=0,1,\ldots,\hat I\}$,
$x\in\cp\cap\La_n$, as described at the beggining of this subsection, let us identify it as 
\begin{equation}\label{eq:ts}
\{\hat s'_i(z);\,i=0,1,\ldots,\hat I\}:=\Xi (\{\hat s_i(x);\,i=0,1,\ldots,\hat I\}), 
\end{equation} 
$z\in\cp'\cap\La'_n$, where for each $z\in\cp'\cap\La'_n$ such that 
$z=\Xi (x)$ for some $x\in\cp\cap\La_n$, we have that $\hat s'_i(z)=\Xi (\hat s_i(x))$, $i=0,1,\ldots,\hat I$. 

It follows from the discussion above that, given $z\in\cp'\cap\La'_n$ and $\hat s'_{i-1}(z)\in\bar\La'_n$ for some
$i\geq1$, if $\ct'_{\hat s'_{i-1}(z),\log n}\cap\cp'\ne\emptyset$, then $\hat s'_{i}(z)$ 
is the point of $\ct'_{\hat s'_{i-1}(z),\log n}\cap\cp'$ whose second coordinate is closest to that of 
$\hat s'_{i-1}(z)$; if on the other hand $\ct'_{\hat s'_{i-1}(z),\log n}\cap\cp'=\emptyset$, then
$\hat s'_{i}(z)$ is the intersection of the vertical straight line segment $(\hat s'_{i-1}(z)_1,\hat s'_{i-1}(z)_2+u)_{u\geq0}$ --- 
where $\hat s'_{i-1}(z)_j$ is the $j$-th coordinate of $\hat s'_{i-1}(z)$, $j=1,2$ --- 
and the top boundary of $\ct'_{\hat s'_{i-1}(z),\log n}$. 
We may characterize $\hat I$ in terms of the $\hat s'_{i}(z)$'s as 
\begin{equation}\label{eq:hi}
\hat I=\inf\{i\geq1:\hat s'_{i}(z)\notin\bar\La'_n\}.
\end{equation}

So, given $\hat\g_x\in\hat\G_n$ for some $x\in\cp\cap\La_n$, in which case $\hat\g_x$ is determined
by $\{\hat s_i(x);\,i=0,1,\ldots,\hat I\}$ and the straight line edges $(\hat s_{i-1}(x),\hat s_i(x))$, 
$i=1,\ldots,\hat I$, the path $\hat\g'_z=\Xi (\hat\g_x)$ is determined by $\{\hat s'_i(z);\,i=0,1,\ldots,\hat I\}$
and the edges $(\hat s'_{i-1}(z),\hat s'_i(z))=\Xi ((\hat s_{i-1}(x),\hat s_i(x)))$, which are {\em not} straight line 
segments.

\paragraph{Modification of $\ct'_{(y,s),l}$.} 

For convenience, for each set $(y,s)\in\cp'\cap\bar\La'_n$ and $l\in[0,\log n]$, we will replace $\ct'_{(y,s),l}$ 
by another set, $\ct''_{(y,s),l}$, such that the area of the symmetric difference is small enough
so that with high probability $\hat\G'_n$ is unchanged by the replacement. We start by observing
that from~(\ref{eq:a}) we have that
\begin{equation}\label{eq:a1}
\frac a{l'}=\frac1n\frac{(1+c^2)(1+\frac sn)}2+O(n^{-2})
\end{equation} 
Substituting that into~(\ref{eq:eta}) we find, after further computation, that
\begin{equation}\label{eq:neta}
n\eta =c\left(1+d/n\right)\frac{1+\frac sn}{1-\frac{l'}n\left(1+\frac sn\right)}\,l'+O((l')^3n^{-2}),
\end{equation} 
where $d=(1+c^2)/[2(1+s/n)]$. So, defining
\begin{equation}\label{eq:ttt}
\ct''_{(y,s),v}=
 \bigcup_{u\in[0,v]}\left[y\pm c\,\frac{1+d/n}{1+s/n}\,u\right]\times
 \left\{s+u\right\},
\end{equation} 
we have that the symmetric difference between the sets $\ct'_{(y,s),\log n}$ and
$\ct''_{(y,s),d'_n\log n}$, say $\Delta_{(y,s)}$, where 
$d'_n=d'_n(s)=\frac{(1+\frac sn)^2}{1-\frac{\log n}n\left(1+\frac sn\right)}$,
has area of order $O((\log n)^4n^{-2})$, which is then also the order of the probability that
we find a point of $\cp'$ in $\Delta_{y,s}$ for any $(y,s)\in\cp'\cap\bar\La'_n$. 
Arguing as is the proof of Lemma~\ref{lm:pp} above, we get that the probability of
finding a point of $\cp'$ in $\Delta_{y,s}$ for some $(y,s)\in\cp'\cap\bar\La'_n$ goes to 0
as $n\to\infty$, and we have that with high probability $\hat\G'_n$ is unchanged if we 
use $\ct''_{(y,s),d'_n\log n}$ instead of $\ct'_{\hat s'_{i-1}(z),\log n}$, $(y,s)\in\cp'\cap\bar\La'_n$,
in its definition. So we will do: for $(y,s)\in\La_n'\cap\cp'$ fixed, 
let us consider $\{\hat s''_i(y,s);\,i=0,1,\ldots,\hat I''\}$, where 
$\hat s''_i(y,s)$ is defined as $\hat s'_i(y,s)$ was in the full paragraph 
below~(\ref{eq:ts}) above, except that we use $\ct''_{(y,s),d'_n\log n}$ instead of 
$\ct'_{\hat s'_{i-1}(z),\log n}$. Then, as just argued, we have that with high probability
for all $(y,s)\in\La_n'\cap\cp'$
\begin{equation}\label{eq:whp} 
\hat I''=\hat I\mbox{ and } 
\hat s''_i(y,s)=\hat s'_i(y,s);\,i=0,1,\ldots,\hat I,
\end{equation}
so we may and will replace $\{\hat s'_i(y,s);\,i=0,1,\ldots,\hat I\}$ 
by $\{\hat s''_i(y,s);\,i=0,1,\ldots,\hat I''\}$, $(y,s)\in\cp'\cap\La'_n$ throughout.

Notice that $\ct''_{(y,s),d'_n\log n}$ is an isosceles triangle similar to the one depicted
in Figure~\ref{fig:rshape}, with height $d'_n\log n$ and base 
$2c\frac{1+\frac dn}{1-\frac{\log n}n\left(1+\frac sn\right)}(1+\frac sn)\log n$ instead.
%$2c(1+d/n)\frac{1+\frac sn}{1-\frac{\log n}n\left(1+\frac sn\right)}\log n$ instead.

Notice also that $\{\hat s''_i(y,s);\,i=0,1,\ldots,\hat I''\}$ (as well as the other {\em paths}
considered so far) may be similarly defined with arbitrary starting point in $\La_n'$ (in this case).

\paragraph{Path increments.}

%We describe now the sequence $\hat s'_i(z);\,i=0,1,\ldots$,

The increments of the paths of $\hat\G'_n$ may be understood as follows. Let us consider the following two independent
families of independent random variables indexed by points of $\bar\La'_n$. Let $\{U_{y,s},\,(y,s)\in\bar\La_n'\}$
be iid uniform in $[-1,1]$ and $\{T_{y,s},\,(y,s)\in\bar\La_n'\}$ be such that
\begin{equation}\label{eq:idf}
\P(T_{y,s}>v)=\exp\left\{-\frac{c_n}{(1+\frac sn)^2}\frac{v^2}{(1+\frac{s+v}n)^2}\right\} 1_{\{v<L_n\}},
\end{equation}
where $c_n:=c(1+\frac dn)$, $L_n:=\frac{(1+\frac sn)^2\log n}{1-\frac{\log n}n(1+\frac sn)}$. And for
$(y,s)\in\bar\La_n'$ let
\begin{equation}\label{eq:xys}
X_{y,s}=\frac{c_n}{1+\frac sn}T_{y,s}1_{\{T_{y,s}<L_n\}}U_{y,s}.
\end{equation}

Let us now define for $(y,s)\in\La_n'\cap\cp'$, $(Y_0,S_0)=(y,s)$, and for $i\geq1$
\begin{equation}\label{eq:ysn}
(Y_i,S_i)=(Y_{i-1},S_{i-1})+(X_{Y_{i-1},S_{i-1}},T_{Y_{i-1},S_{i-1}}),
\end{equation}
and let
\begin{equation}\label{eq:i}
\ci=\inf\{i\geq1:\,(Y_i,S_i)\notin\bar\La'_n\}.
\end{equation}

\begin{lem}\label{lm:idd}
Given  $(x,t)\in\La_n'$, we have that
\begin{equation}\label{eq:idd}
\{\hat s''_i(x,t);\,i=0,1,\ldots,\hat I''\}\stackrel{d}{=}\{(Y_i,S_i);\,i=0,1,\ldots,\ci\},
\end{equation}
with $(Y_0,S_0)=(x,t)$, where $\stackrel{d}{=}$ means identity in distribution.
\end{lem}
\noindent{\bf Proof} 

Follows from elementary properties of Poisson point processes. 

Notice that given $\hat s''_{i-1}(x,t)=(y,s)$ and $0\leq v<L_n$, the event 
$\{\hat s''_{i}(x,t)_2-\hat s''_{i-1}(x,t)_2>v\}$ corresponds to 
$\{\ct''_{(y,s),v}\cap\cp'=\emptyset\}$; the probability of the latter
event is then the exponential of minus the integral over $\ct''_{(y,s),v}$
against the intensity measure of $\cp'$, given in~(\ref{eq:int}); a straightforward
computation of that integral gives the absolute value of the expression within 
braces in~(\ref{eq:idf}) in that case. The event 
$\{\hat s''_{i}(x,t)_2-\hat s''_{i-1}(x,t)_2=L_n\}$ corresponds to 
$\{\ct''_{(y,s),L_n}\cap\cp'=\emptyset\}$ and the expression in~(\ref{eq:idf})
follows in this last case. Now given $\hat s''_{i-1}(x,t)=(y,s)$ and 
$\hat s''_{i}(x,t)_2-\hat s''_{i-1}(x,t)_2=v$, $0\leq v\leq L_n$, it follows from 
elementary properties of Poisson point processes that 
$\hat s''_{i}(x,t)_1-\hat s''_{i-1}(x,t)_1$ is either uniformly distributed 
on the top side of $\ct''_{(y,s),v}$, if $v<L_n$, or it vanishes, if $v=L_n$, 
giving rise to the distribution of $X_{y,s}$ in both cases.
$\square$

\medskip

We are now in the position to drop a modification operated in the construction 
of our paths, namely the one involving lateral excursions of 
$\hat s''_i(y,s);\,i=0,1,\ldots,$ outside the left or right lateral boundaries 
of $\bar\La'_n$ till it reaches the top of $\bar\La'_n$. We show next that
with high probability these excursions do not happen for any $(y,s)\in\La_n'\cap\cp'$.

\begin{lem}\label{lm:lat}
Given $(y,s)\in\La_n'\cap\cp'$, let 
\begin{equation}\label{eq:J}
J=\inf\{i\geq1:\,\hat s''_i(y,s)_2>\tau n\}.
\end{equation}
Then with high probability
$|\hat s''_i(y,s)_1|\leq n^{1-a}$ (simultaneously) for all $i=1,\ldots,J$
and all $(y,s)\in\La_n'\cap\cp'$.
\end{lem}

\noindent{\bf Proof} 

It follows from~(\ref{eq:idf}) and Lemma~\ref{lm:idd} that $\{T_{y,s};\,(y,s)\in\bar\La'_n\}$ 
dominates stochastically an iid family $\{V_{y,s};\,(y,s)\in\bar\La'_n\}$ such that
\begin{equation}\label{eq:idf1}
\P(V_{y,s}>v)=e^{-Cv^2} 1_{\{v<L_n\}},
\end{equation}
for some constant $C$. Standard large deviation estimates then yield the existence of positive constants
$C', C''$ such that 
\begin{equation}\label{eq:ld}
\P(J>C'n)\leq e^{-C''n}.
\end{equation}

Conditioning on $\hat s''_i(y,s)_2$, $i\geq1$, such that $J\leq C'n$, and using~(\ref{eq:xys}) 
and Lemma~\ref{lm:idd}, we have that
\begin{equation}\label{eq:lat1}
 \P'\left(\max_{1\leq i\leq J}(\hat s''_i(y,s)_1-y)>n^{1-b}\right)
 \leq\sum_{1\leq i\leq J}\P'(\hat s''_i(y,s)_1-y>n^{1-b}),
\end{equation}
where $\P'$ is the appropriate conditional distribution. Resorting standard large deviation
estimates, the latter probability is bounded above by
\begin{eqnarray}\label{eq:lat2}
e^{-\l n^{1-b}}\prod_{i=1}^J\E(e^{\l T'_iU}),
\end{eqnarray}
where $\l>0$ is a positive number to be chosen, and 
$T_i'=\frac{c_n}{1+{S''_{i-1}}/{n}}T''_i1_{\{T''_i<L_n\}}$, with $S''_i=\hat s''_i(y,s)_2$ e 
$T''_i=S''_i-S''_{i-1}$.
Now the expected value in~(\ref{eq:lat2}) is bounded above by
\begin{eqnarray}\label{eq:lat3}
\E(e^{c_n\l T''_iU})=\frac{e^{\l T''_i}-e^{-\l_n T''_i}}{2\l_n T''_i}\leq 1+\l_n^2 (T''_i)^2,
\end{eqnarray}
where $\l_n=1/\sqrt n$, and we have made the choice $\l=1/(c_n\sqrt n)$, and we have made use 
of the fact that $T''_i\leq L_n$, and thus $\l_n T''_i=o(1)$ in the latter inequality of~(\ref{eq:lat3}).
The product in~(\ref{eq:lat1}) is thus bounded above by the exponential of
\begin{eqnarray}\label{eq:lat4}
\l_n^2 \sum_{i=1}^J(T''_i)^2\leq2(1+\tau)^2\log n\,\l_n^2\sum_{i=1}^JT''_i
%\leq2\l_n^2L_n(\tau n+(1+\tau)^2\log n,
\end{eqnarray}
and we estimate the latter sum as
\begin{eqnarray}\label{eq:lat5}
\sum_{i=1}^{J-1}T''_i+T''_J\leq\tau n+2(1+\tau)^2\log n\leq 2\tau n
\end{eqnarray}
for all large enough $n$.

Replacing now our estimates in~(\ref{eq:lat2}), we get the following upper bound for the sum 
%conditional probability
on the left hand side of~(\ref{eq:lat1}) for all large $n$, uniformly in the conditioning variables (within the
prescribed conditions).
\begin{eqnarray}\label{eq:lat6}
 C'n\exp\{-c_n^{-1}n^{\frac12-b}+4\tau(1+\tau)^2\log n\}\leq\exp\left\{-\frac1{2c}n^{\frac12-b}\right\}
\end{eqnarray}
Using~(\ref{eq:ld}), we find that, for all large enough $n$, twice that %a similar 
bound follows for the expression that we obtain when we replace the conditional $\P'$ by 
the unconditional $\P$ in the left hand side of~(\ref{eq:lat1}). The same argument then provides the same bound
when we replace $(\hat s''_i(y,s)_1-y)$ by $-(\hat s''_i(y,s)_1-y)$ in that probability, and we thus get 
twice the bound for $|\hat s''_i(y,s)_1-y|$. 

Finally, repeating an argument already used in the proof of Lemma~\ref{lm:id} concerning the cardinality of
Poisson points within a region, this time the cardinality of $\cp'\cap\La'_n$, we have that outside 
an event of vanishing probability as $n\to\infty$, it is less than constant times $n^{2-b}$. 
Using that and the above estimates, we then have that 
\begin{equation}\label{eq:lat7}
 \P\left(\max_{(y,s)\in\cp'\cap\La'_n}\max_{1\leq i\leq J}|\hat s''_i(y,s)_1-y|>n^{1-b}\right)
 \leq 4C''' n^{2-b}e^{-\frac{n^{-b+1/2}}{2c}}+\P(|\cp'\cap\La'_n|>C'''n^{2-a})
\end{equation}
vanishes as $n\to\infty$, provided the constant $C'''$ is large enough, and since $b<1/2$.

Now outside the event in the probability in the left hand side of~(\ref{eq:lat7}), we have that
\begin{equation}\label{eq:lat8}
\max_{(y,s)\in\cp'\cap\La'_n}\max_{1\leq i\leq J}|\hat s''_i(y,s)_1|
\leq n^{1-b}+\max_{(y,s)\in\cp'\cap\La'_n}|y|\leq 2n^{1-b}\leq n^{1-a}
\end{equation}
for all large enough $n$, and we are done.
$\square$

\paragraph{Variant of $\hat\G'_n$.}

Summing the above up, we may and will replace the paths of $\hat\G'_n$ determined by 
$\{\hat s''_i(y,s);\,i=0,1,\ldots,\hat I''\}$ by those determined by 
$\{\hat s''_i(y,s);\,i=0,1,\ldots,J\}$, $(y,s)\in\La_n'\cap\cp'$. Indeed, for 
$(y,s)\in\La_n'\cap\cp'$, let $\hat\g''_{y,s}$ be the path determined by 
$\{\hat s''_i(y,s);\,i=0,1,\ldots,J\}$ as suggested above: $\hat\g''_{y,s}$
starts at $(y,s)$ and runs through the {\em edges} $(\hat s''_{i-1}(y,s),\hat s''_i(y,s))$,
$i=1,\ldots,J-1$ and the last edge $(\hat s''_{J-1}(y,s),\hat s^*_J(y,s))$, 
where $\hat s^*_J(y,s)$ is the point of the top side of $\bar\La'_n$ intersected by
the edge $(\hat s''_{J-1}(y,s),\hat s''_J(y,s))$. The edge $(\hat s''_{i-1}(y,s),\hat s''_i(y,s))$
is given by $(\hat s'_{i-1}(y,s),\hat s'_i(y,s))$, whenever the respective endpoints of the pair of
edges coincide --- in this case, it is not linear, as remarked above; see discussion at the end of 
paragraph below~(\ref{eq:hi}) ---; otherwise (an event such that the union of all such events as the
start point varies over $\La_n'\cap\cp'$ has vanishing probability as $n\to\infty$), 
it is the linear interpolation of its endpoints. 

Let then
\begin{equation}\label{eq:hgpp}
\hat\G''_n=\{\hat\g''_{y,s};\,(y,s)\in\La_n'\cap\cp'\}.
\end{equation}

From the above discussion, we may replace $\hat\G'_n$ by $\hat\G''_n$ in proving Proposition~\ref{aux}
and Lemma~\ref{lm:dist}.
Namely, it is enough to prove the following results.

Let 
\begin{equation}\label{eq:reschgpp}
 \tilde\G''_n=\{D(\hat\g'');\,\hat\g''\in\hat\G''_n\}
\end{equation}
be the collection of diffusively rescaled paths of $\hat\G''_n$ --- see~(\ref{eq:tt}).

\begin{prop}\label{vaux}
As $n\to\infty$
 \begin{equation}
  \label{eq:vaux}
  \tilde\G''_n\Rightarrow\cw^\tau_0.
 \end{equation}
 %where $\cw^\tau_0$ is the restricted Brownian web given in~(\ref{rbw}) above, and ``$\Rightarrow$'' 
 %stands for convergence in distribution.
\end{prop}

We will prove this result in the remaining sections of this paper.

\begin{lem}\label{lm:vdist}
Let $\a'\in(0,\a)$ be given. We have
\begin{equation}
\label{eq:vdist}
d_{\ch_{-1}^{-\a,\a'}}(\tilde\G_n,\psi''(\tilde\G''_n))\to0
\end{equation} 
with high probability as $n\to\infty$.
\end{lem}

We will prove this result in the next subsection.

%%%%%%%%%%%%%%%%%%%%%%%%%%%%%%%%%%%%%%%%%%%%%%%%%%%%%%%%%%%%%%%%%%%%%%%%%%%%%%%%%%%
%%%%%%%%%%%%%%%%%%%%%%%%%%%%% PROOF %%%%%%%%%%%%%%%%%%%%%%%%%%%%%%%%%%%%%%%%%%%%%%
%%%%%%%%%%%%%%%%%%%%%%%%%%%%%%%%%%%%%%%%%%%%%%%%%%%%%%%%%%%%%%%%%%%%%%%%%%%%%%%%%%

\subsection{Proof of Lemma~\ref{lm:vdist}}\label{pflm}

%\setcounter{equation}{0}

%\subsection{Preliminaries}\label{pre}

%\paragraph{}

Using the notation of Section~\ref{mod}, we will show that with high probability, given a polygonal path $\g$
determined by $\{s_i(x),\,i=0,1,\ldots, I'\}$, $x\in\cp\cap\La_n$, then letting $\tilde\g=D(\g)$ be the diffusive
rescaling of $\g$, and $\mathring\g=\psi\circ D\circ\Xi(\g)$ be the image under $\psi$ of the diffusive rescaling 
of the image under $\Xi$ of $\g$, we have that 
\begin{enumerate}
 \item as subsets of the plane, $\tilde\g$ and $\mr\g$ belong to $[-n^{1-a},n^{1-a}]\times[-1,-\a']$; and
 \item $d(\tilde\g,\mr\g)\to0$ as $n\to\infty$ uniformly over $x\in\cp\cap\La_n$,
\end{enumerate}
where $d$ is the distance in $\Pi_{-1}^{-\a,-\a'}$ defined
in~(\ref{metric}) above.

The first point follows from our arguments above, in particular those in the proofs of Lemmas~\ref{lm:id} 
and~\ref{lm:lat}. We will assume below that the properties, which from the latter lemmas hold with high probability, 
are in force.

As for the second point, let us start by considering the vertices of the polygonal line forming $\tilde\g$, namely
$D(s_i(x)),\,i=0,1,\ldots, I'$. Let $w=(w_1,w_2)$ be one of those points. Let us write 
$w=D(re^{-\frac\pi2+\s})=\le(\frac{r\sin\s}{\sqrt n},-\frac{r\cos\s}n\ri)$, where $re^{-\frac\pi2+\s}=s_i(x)$
for some $i=0,1,\ldots, I'$. Then, $\psi\circ D\circ\Xi(re^{-\frac\pi2+\s})
=\mr w=(\mr w_1,\mr w_2)=\le(\frac{r\s}{\sqrt n},-\frac{r\s}n\ri)$ is one of the
{\em vertices} of $\mr\g$. Since $r\leq n$ and $|\s|\leq n^{-a}$, we have that 
\begin{equation}\label{eq:mrw}
|w_1-\mr w_1|=O(n^{-3a+1/2}),\quad|w_2-\mr w_2|=O(n^{-2a}) 
\end{equation}
uniformly over $x\in\cp\cap\La_n$. It immediately follows that the absolute value 
of the difference between the starting times of $\tilde\g$ and $\mr\g$ is $o(1)$ uniformly over $x\in\cp\cap\La_n$,
and the same holds for the ending times. It is enough then to have the same estimate for 
$\sup_{-1\leq s\leq\a'}|\tilde\g^\star(s)-\mr\g^\star(s)|$, where $\tilde\g^\star$ and $\mr\g^\star$ are
$\tilde\g$ and $\mr\g$ respectively continued to $[-1,\a']$, as defined in the paragraph of~(\ref{metric}) above.

In order to acomplish that, we first notice that, with the notation of the previous paragraph, $w_2\geq\mr w_2$.
%Let us consider $\mr\g^\star(w_2)$. If 
Suppose $w=D(s_{I'}(x))$, and take $s\geq w_2$. Then $\tilde\g^\star(s)=w_1$, $\mr\g^\star(s)=\mr w_1$, and thus
$|\tilde\g^\star(w_2)-\mr\g^\star(w_2)|=|w_1-\mr w_1|$, wich as established in the previous paragraph,
equals $O(n^{-3a+1/2})$ uniformly over $x\in\cp\cap\La_n$. 

To bound $|\tilde\g^\star(s)-\mr\g^\star(s)|$ when $s=w_2$ for $i<I'$, we first claim that
the portion of $\mr\g^\star$ above $\mr w_2$ and up to $w_2$ is contained in an isosceles triangle
like the one depicted in Figure~\ref{fig:rshape} with lower vertex at $\mr w$, internal angle at $\mr w$
whose tangent equals $c'\sqrt n$, $c'$ a constant, and height $w_2-\mr w_2$. This follows from a straightforward
analysis of $\psi\circ D(\ct''_{\mr w,w_2-\mr w_2})$; see~(\ref{eq:xt}), (\ref{eq:eta}) and~(\ref{eq:neta}).
We then have by~(\ref{eq:mrw}) that 
$|\tilde\g^\star(w_2)-\mr\g^\star(w_2)|\leq|w_1-\mr w_1|+c''{\sqrt n}|w_2-\mr w_2|\leq c'''n^{-2a+1/2}$
for some constants $c'',c'''$.

Finally, for $s$ between succesive $w_2$'s (say, $w_2$ and $w_2'$, with $w_2<w_2'$), we may similarly as
in the previous paragraph argue that the portion of $\mr\g^\star$ above $\mr w_2$ and up to $s$ is contained in a
triangle like the one of the previous paragraph, except that with height larger by $s-\mr w_2\leq$ const $\log n/n$
than the one of the previous paragraph, and the same bound (with a larger $c'''$) follows for 
$|\tilde\g^\star(s)-\mr\g^\star(s)|$. The argument is thus complete, once we recall that $a>1/4$. $\square$

%%%%%%%%%%%%%%%%%%%%%%%%%%%%%%%%%%%%%%%%%%%%%%%%%%%%%%%%%%%%%%%%%%%%%%%%%%%%%%%%%
%%%%%%%%%%%%%%%%%%%%%%%%%%%%% PROOF %%%%%%%%%%%%%%%%%%%%%%%%%%%%%%%%%%%%%%%%%%%%%%
%%%%%%%%%%%%%%%%%%%%%%%%%%%%%%%%%%%%%%%%%%%%%%%%%%%%%%%%%%%%%%%%%%%%%%%%%%%%%%%%%%

\section{Proof of Proposition~\ref{vaux}}\label{pfpr}

\setcounter{equation}{0}

%\subsection{Preliminaries}\label{pre}

%\paragraph{}

For the convenience of being able to deal with paths starting from any point of $\La_n'$, 
we will consider the following {\em completion} of $\hat\G''_n$.
\begin{equation}\label{eq:hgppp}
\hat\G'''_n=\{\hat\g'''_{x};\,x\in\La_n'\},
\end{equation}
where $\hat\g'''_{x}$ is a path determined by $\{\hat s''_i(x);\,i=0,1,\ldots,J\}$, 
using $\ct''_{x,d'_n\log n}$, as for $\hat\g''_{x}$, except that now it may start from
any $x\in\La_n'$. We will stipulate that the edges of $\hat\g'''_{x}$ are straight line 
segments until a point of $\cp'$, say $x'$, is hit, after which the edges are identical with 
those of $\hat\g'''_{x'}$.

\paragraph{Claim.}
We now claim that with high probability $\hat s'''_1(x)\in\cp'$
%the end of the first edge of $\hat\g'''_{y,s}$ belongs to $\cp'$ 
for all $x\in\La_n'$. Indeed, let us fix $x\in\La_n'$
and take $x_n$ as a closest point to $x$ of $\Z^2\cap\ct''_{x,d'_n\log n}$. We note that
for all  $n$ large enough $\ct''_{x_n,\frac12d'_n\log n}\subset\ct''_{x,d'_n\log n}$
for all $x\in\La_n'$. The claim will be established once we show that the following probability
\begin{equation}
 \P(\cup_{x\in\Z^2\cap\La_n'}\{\cp'\cap\ct''_{x,d'_n\log n}=\emptyset\})
\end{equation}
vanishes as $n\to\infty$. But this probability is bounded above by constant times
$n^{2-a}$ times $\exp\{-\mbox{area of $\ct'_{x,\frac12\log n}$}\}$. Since the latter area
is bounded below by constant times $(\log n)^2$, the claim follows.

We may then suppose that $\hat s'''_1(x)\in\cp'$ for all $x\in\La_n'$.
Notice that if $x\notin\cp'$, then  $\hat\g'''_{x}$ coincides
with $\hat\g''_{\hat s''_1(x)}$ from its second vertex -- given by $\hat s''_1(y,s)$ -- on.
Otherwise $\hat\g'''_{x}=\hat\g''_{x}$.

Let 
\begin{equation}\label{eq:reschgppp}
 \tilde\G'''_n=\{D(\hat\g''');\,\hat\g'''\in\hat\G'''_n\}.
\end{equation}
Below, for $x\in\R\times[0,\tau]$, we will write $\tilde\g'''_x$ for 
$D(\hat\g'''_{D^{-1}(x)})=D(\hat\g'''_{(\sqrt{n}x_1,nx_2)})$.

We will establish the following result.
\begin{prop}\label{vauxa}
As $n\to\infty$
 \begin{equation}
  \label{eq:vauxa}
  \tilde\G'''_n\Rightarrow\cw^\tau_0.
 \end{equation}
 %where $\cw^\tau_0$ is the restricted Brownian web given in~(\ref{rbw}) above, and ``$\Rightarrow$'' 
 %stands for convergence in distribution.
\end{prop}

This immediately implies the claim of Proposition~\ref{vaux}, once we show the following lemma.

\begin{lem}\label{lm:vdista}
We have that
\begin{equation}
\label{eq:vdista}
d_{\ch_{0}^{\tau,\tau}}(\tilde\G'''_n,\tilde\G''_n)\to0
\end{equation} 
with high probability.
\end{lem}

\noindent{\bf Proof of Lemma~\ref{lm:vdista}} 

Given $x=(y,s)\in\La_n'$ and $\hat\g'''_{x}\in\hat\G'''_n$, if $s\leq \tau n-d_n'\log n$, then,
as noticed above, we have that (with high probability simultaneously for all $x\in\La_n'$)
$\hat\g'''_{x}$ coincides with $\hat\g''_{x'}\in\hat\G''_n$ from its second edge $x'$ on.
We thus have that the distance of starting times of $\tilde\g'''_{x}$ and $\tilde\g''_{x'}$
is bounded above by constant times $n^{-1}\log n$, and the (uniform) distance between
the continuations (starred versions) of $\tilde\g'''_{x}$ and $\tilde\g''_{x'}$, 
which is then given by $n^{-1/2}|x_1-x_1'|$, which is bounded above by constant times 
$n^{-1/2}\log n$.

Let us now examine the case where $y>\tau n-d_n'\log n$. Let us introduce 
$\hat\ct''_{x,d'_n\log n}$, the triangle obtained by reflecting $\ct''_{x,d'_n\log n}$
on the horizontal axis through $x$. Then again we may argue that 
with high probability $\ct''_{x,d'_n\log n}\cap\cp'\ne\emptyset$ for all $x\in\La_n'$.
Let $x''$ be the point of $\ct''_{x,d'_n\log n}\cap\cp'$ closest to $x$, and let us 
consider $\tilde\g'''_{x}$ and $\tilde\g''_{x''}$. The distance between their starting points
is thus bounded above by $n^{-1}\log n$. Both paths are contained in isosceles triangles 
shaped like $\ct''_{x,d'_n\log n}$, but with the tangent to the internal angle at the bottom vertices
($x$ and $x'$, respeectively) a constant times $\sqrt n$, and height a constant times $\log n/n$.
We conclude that the uniform distance between the continuations of $\tilde\g'''_{x}$ and $\tilde\g''_{x''}$
is bounded above by constant times $n^{-1/2}\log n$.

We conclude from the above paragraphs for every path $\tilde\g'''\in\tilde\G'''_n$, we may find a
path $\tilde\g''\in\tilde\G''_n$ such that $d(\tilde\g''',\tilde\g'')$ is bounded above by (a uniform)
constant times $n^{-1/2}\log n$, and we are done
$\square$

\medskip

In order to prove Proposition~\ref{vauxa}, we will verify criteria of~\cite{finr1}, which were devised
for the  case of convergence to the ordinary Brownian web, adapted in an obvious way for convergence 
to the restricted Brownian web. 
%(That the adaptation leads to valid criteria follows since $\ch_0^{\tau,\tau}$ 
%is the image under a continuous map of a closed subset of $\bar\ch$, the Brownian web sample space, and the metric on 
%$\ch_0^{\tau,\tau}$ we are considering is equivalent to that induced by the metric on $\bar\ch$ via the 
%mentioned continuous map. See discussion in Section~\ref{mod}.)
The adapted criteria are as follows.

\begin{itemize}
\item[($I$)] Let $\mathcal{D}$ be a countable dense deterministic set in $\R\times(0,\tau)$. For any %deterministic 
$x_1,...,x_m \in \mathcal{D}$, $\tilde{\gamma}'''_{x_1}, ... , \tilde{\gamma}'''_{x_m} \in \tilde{\G}'''_n$ 
converge in distribution as $n \rightarrow \infty$ to coalescing Brownian motions, with the same diffusion 
constant $\o$, starting from $x_1,...,x_m$, and ending at time $\tau$. 
\item[($B_1^\prime$)] For a system $\mathcal{V}$ of space-time trajectories in $\mathbb{R}^2$, let 
$\eta_{\mathcal{V}}(t_0,t;a,b)$, $a<b$, denote the random variable that counts the number of distinct points 
in $\mathbb{R}\times \{t_0+t\}$ that are touched  by paths in $\mathcal{V}$ which also touch some point in 
$[a,b]\times \{t_0\}$. Then, for every $0<\beta <\tau$,
$$
\limsup_{\epsilon \rightarrow 0+} \limsup_{n \rightarrow +\infty} \sup_{ t > \beta } 
\sup_{ (a,t_0) \in \mathbb{R}\times[0,\tau-\b] } 
\textrm{P} ( \eta_{\tilde{\G}'''_n}(t_0,t;a,a+\epsilon) \ge 2) = 0 \, .
$$
\item[($B_2^\prime$)] For every $\beta > 0$,
$$
\limsup_{\epsilon \rightarrow 0+} \frac{1}{\epsilon} \limsup_{n \rightarrow +\infty}
\sup_{ t > \beta } \sup_{ (a,t_0) \in \mathbb{R}\times[0,\tau-\b] } 
\textrm{P} ( \eta_{\tilde{\G}'''_n} (t_0,t;a,a+\epsilon) \ge 3) = 0 \, .
$$
\end{itemize}

That these criteria are valid (as criteria) can be argued with a straightforward and obvious adaptation of the 
arguments for the proof of Theorem 2.2 in~\cite{finr1}.

\medskip

It turns out that we do not know how to verify $B_2^\prime$ for the present model. Resorting to FKG-type inequalities,
a strategy that has been successful in a few cases, does not seem to be a way.
We resort instead to an alternative criterium to $B_2^\prime$, as adopted in~\cite{nrs} in the study
of a coalescing system in the full space-time plane. We present it next, again in a suitable adaptation for the case
of $\R\times[0,\tau]$. 
%a condition is derived to replace $(B^\prime_2)$ in the proof of convergence of diffusively rescaled 
%non-nearest neighbor coalescing random walks to the Brownian web. Their condition is based on duality and can be expressed as
\begin{itemize}
\item[($E$)] For a system $\mathcal{V}$ of space-time trajectories in $\mathbb{R}\times[0,\tau]$, let 
$\hat{\eta}_\mathcal{V} (t_0,t;a,b)$, $a<b$, $t_0\in[0,\tau)$, $t\in(0,\tau-t_0]$, be the number distinct 
points in $(a,b)\times \{t_0+t\}$ that are touched by a 
path that also touches $\mathbb{R} \times \{t_0\}$. Then for any subsequential limit $\mathcal{X}$ of $\tilde{\G}'''_n$ we have that
\begin{equation}\label{cond_e}
 \textrm{E}[ \hat{\eta}_\mathcal{X} (t_0,t;a,b) ] \leq \frac{b-a}{\sqrt{\pi \o^{-1} t}} \, . 
\end{equation} 

(In \cite{nrs} the diffusion coefficient of the Brownian web is 1, so the factor  $\o^{-1}$ does not appear (explictly). 
It is of course straightforward to go from that condition to~(\ref{cond_e}) by rescaling time suitably.)
\end{itemize}

Again, the legitimacy of ($E$) as a criterium may be argued entirely as in~\cite{nrs}, with obvious adaptations. 

\medskip

The verification of condition $E$ requires tightness of $\hat{\G}'''_n$, but this follows from condition $I$ due to the noncrossing 
property (see Proposition B2 in \cite{finr1}). Futhermore, once we have proper estimates on the distributions of the coalescence times 
and condition $I$, the verification of conditions $B^\prime_1$ and $E$ follows from adaptations of arguments presented 
in~\cite{nrs},~\cite{s} and~\cite{cfd}, as we will see below. 
The estimate we need on the distribution of the coalescence time of two random trajectories of $\hat{\G}'''_n$ is 
that the probability that they coalesce after time $t$ is of order $1/\sqrt{t}$. 
This is the content of Proposition~\ref{prop:coaltime1} below.

%%%%%%%%%%%%%%%%%%%%%%%%%%%%%%%%%%%%%%%%%%%%%%%%%%%%%%%%%%%%%%%%%%%%%%%%%%%%%%%%%
%%%%%%%%%%%%%%%%%%%%%%%%%%%%% Glauco %%%%%%%%%%%%%%%%%%%%%%%%%%%%%%%%%%%%%%%%%%%%%%
%%%%%%%%%%%%%%%%%%%%%%%%%%%%%%%%%%%%%%%%%%%%%%%%%%%%%%%%%%%%%%%%%%%%%%%%%%%%%%%%%%

\section{Coalescence time}
\label{sec:coaltime}

\setcounter{equation}{0}

One important step to establish convergence in distribution of a family of random coalescing paths to the BW is to prove that 
the tail of the distribution of the coalescence time between two such paths decays as $O(1/\sqrt{t})$. In this section, we want 
to avoid problems with the fact that the paths of $\hat{\Gamma}_n'''$, $n\geq1$, are defined on bounded time intervals. 
(Recall that $\hat{\Gamma}_n'''$ is defined in the interval $[0,\tau n]$, and so is the Poisson point process $\mathcal{P}^\prime_n$). 
So we begin by extending $\mathcal{P}^\prime_n$ to a Poisson point process on the upper half plane $\H:=\mathbb{R}\times[0,\infty)$, 
and then extend $\hat{\Gamma}_n'''$ accordingly, as follows.

Let $\mathcal{P}^{\star}_n$ be a Poisson point process on $\H$ with intensity measure
\begin{equation}
\label{intensityPstar}
\left(\frac{1_{(0,\tau n]}(s)}{(1+s/n)^3} +
\frac{1_{(\tau n, + \infty]}(s)}{(1+\tau)^3} \right) dy ds \, . 
\end{equation}
So  $\mathcal{P}^{\star}_n$ can be considered the union of $\mathcal{P}^{\prime}_n$ with an independent homogeneous Poisson point 
processes in  $\mathbb{R} \times (\tau n,+\infty)$ with intensity $(1+\tau)^{-3}$.

The paths in $\hat{\Gamma}_n'''$ are restricted to $\Lambda_n^\prime$. We will also drop this restriction. 
But we still need to define $\mathcal{T}_{(y,s)}^{\prime \prime}$ for $(y,s) \in \mathbb{R} \times [\tau n, + \infty]$. 
We take a %the only possible 
natural definition:
\begin{equation}
 \label{t_alt}
\mathcal{T}_{(y,s), d_n^\prime\log n}^{\prime \prime} = 
\mathcal{T}_{(y,s)}^{\prime \prime} = 
\bigcup_{u \in [0,d_n^\prime\log n]} \left[ y \pm c \, \frac{1 + d/n}{1 + \tau} \, u \right] \times
\{s+u\} \ , \quad \textrm{if } s \ge \tau n \, .
\end{equation}

Now let us define the system of random paths that we are going to consider throughout the rest of this section. 
We define $\{\hat{s}^\star_i(x) : i \ge 1 \}$ using $\mathcal{T}_{(y,s), d_n^\prime log(n)}^{\prime \prime}$ 
analogously to $\{ \hat{s}^{\prime \prime}_j : j \ge 1 \}$, but without restricting to $\mathbb{R}\times[0,\tau n]$. 
Let
$$
\hat{\Gamma}^{\star}_n = \{ \hat{\gamma}^\star_x : x \in \H \} \, ,
$$
where $\hat{\gamma}_x^\star$ is a path determined by $x$ and the transition points $\{\hat{s}^\star_i(x) : i \ge 1 \}$ 
as in $\hat{\gamma}_x^{\prime \prime \prime}$ except that now there is no truncation of paths. 

\medskip

For some given $t_0\geq0$, let $ X^0_s = \gamma^{\star}_{(0,t_0)}(t_0+s)$, $s\geq0$, and $ X^m_s = \gamma^{\star}_{(m,t_0)}( t_0+s)$, 
$s\geq0$, be two paths in a given realization of $\hat{\Gamma}^{\star}_n$ starting at time $t_0$ respectively in $0$ and 
$m$, where $m$ is an arbitrary deterministic positive real number. %, such that $bar X^0_{0} = 0$ and $bar X^m_{0} = m$. 
Denote
$$
\nu_m = \inf\{ t > 0 : X^m_{t} - X^0_t = 0 \} \, .
$$ 
The aim of this section is to prove the following proposition:

\medskip

\begin{prop} 
\label{prop:coaltime1}
There exists a constant $C>0$, such that 
\begin{equation}
\label{eqprop:coaltime}
\textrm{P} (\nu_m > t) \leq \frac{C \, m}{\sqrt{t}} \, , \, \textrm{ for every } t>0.
\end{equation}
\end{prop}

\medskip

For the sake of simplifying the notation, we supress $n$ from the notation. However, it is important to point out that the estimates are uniform in $n$. 
In particular, the constant in the statement of Proposition \ref{prop:coaltime1} does not depend on $n$.

\bigskip

\noindent \textbf{Proof}

We start by introducing a jump version of $X^0_t$ and $X^m_{t}$. For $j=0,m$, 
let $\bar X^j_t= X^j_{\hat s^\star_i(j,t_0)}$ whenever $t\in[\hat s^\star_i(j,t_0),\hat s^\star_{i+1}(j,t_0))$ for 
some $i\geq0$. Note that $\nu_m = \inf\{ t > 0 : \bar X^m_{t} - \bar X^0_t = 0 \} \, .$

Define 
$$
S^0_1 = d_n^\prime \log n \wedge \inf\{t> 0 : (\bar X^0_t,t_0+t) \in \mathcal{P}^{\star}_n \} 
, \
S^m_1 = d_n^\prime \log n \wedge \inf\{t> 0 : (\bar X^m_t,t_0+t) \in \mathcal{P}^{\star}_n \} 
$$
and, by induction, for $n \ge 2$, 
\begin{eqnarray*}
 S^0_i \= d_n^\prime \log n \wedge \inf\{t- S^0_{i-1} : t> S^0_{i-1},\,\bar X^0_t \in \mathcal{P}^{\star}_n \}\,, 
\\
S^m_i \= d_n^\prime \log n \wedge \inf\{t- S^m_{i-1} : t> S^m_{i-1}\,,\bar X^m_t \in \mathcal{P}^{\star}_n \} \, .
\end{eqnarray*}
Note that $\sum_{j=1}^iS^0_j$ (resp. $\sum_{j=1}^iS^m_i$) is equal to the second coordinate of the transition point 
$\hat{s}^\star_i((0,t_0))$ (resp. $\hat{s}^\star_i((m,t_0))$).

We have that, for $i=0,m$, $(S^j_i)_{i\ge 0}$ is a sequence of independent random variables. 
Note that they are not identically distributed due to the non-homogeneity of $\mathcal{P}_n^{\star}$. 
Let $(N^j_t)_{t \ge 0}$ be the counting process associated to $(S^j_i)_{i\ge 0}$, $j=0$, $m$, i.e., 
\begin{equation}%\label{pp}
\label{Njt}
N^j_t = \sup \Big\{ k \ge 0 : \sum_{j=0}^k S^j_i \leq t \Big\}  \, .
\end{equation}
Thus $N^j_t$ represents the number of transition points in the time interval $[0,t]$ of the path $(\bar X^j_t)_{t\ge 0}$, $j=0,m$. 
Denote by $\tilde{N}_t$ the total number of transition points in the time interval $[0,t]$ of both paths $(\bar X^0_t)_{t\ge 0}$ 
and $(\bar X^m_t)_{t\ge 0}$. Let $\tilde{S}_1$ be the random time of the first transition point of $\tilde{N}_t$, and $\tilde{S}_i$, 
$i\ge 2$, be the random time %between the (i-1)-st and 
of the i-th jump of $\tilde{N}_t$. 

%\medskip

Put $Z^m_0 = m$ and $Z^m_j = \bar X^m_{\tilde{S}_j} - \bar X^0_{\tilde{S}_{j}}$, $j \ge 1$. Note that
$$
\nu_m = \inf \big\{ t > t_0 : Z^m_{\tilde{N}_t} = 0 \big\} \, .
$$

%\bigskip

Note that $(Z^m_j)_{j \ge 0}$ is a martingale with bounded increments. 
By the Skorohod Embedding Theorem, see \cite{monroe}, there exist a standard Brownian motion $(B(s))_{s\ge 0}$ adapted to a 
certain filtration $(\mathcal{G}_s)_{s\ge 0}$ and stopping times $T_1$, $T_2$, ..., such that 
$Z^m_0 = B(0) = m$ and $Z^m_j = B(T_j)$, for $j \ge 1$. Furthermore, denoting $T_0 = 0$, the stopping times $T_1$, $T_2$, ..., 
have the following representation:
$$
T_n = \inf \big\{ s \ge T_{n-1} : B(s) - B(T_{n-1}) \notin \big(U_n,V_n) \big) \big\} \, ,
$$
where $\{(U_n,V_n): n \ge 1 \}$ is a family of random vectors taking values in $\{(0,0)\} \cup (-\infty,0) \times (0,\infty)$. 
In our case, contrary to the cases analyzed in \cite{cfd,cv,v}, where a similar approach is taken,
the random vectors $(U_n,V_n)_{n \ge 1}$ are not independent. 
%Note also that we have started with $Z^m_0 = B(0)$, although $(Z^m_j)_{j \ge 0}$ is not necessarily a martingale. 
%The martingale property is not required for the Skorohod Embedding of discrete time processes, see Lemma 4 in \cite{monroe}. 

We have the inequality
\begin{equation} \label{ineqct}
\textrm{P} ( \nu_m > t) \leq \textrm{P}\big( \nu_m > t, T_{\tilde{N}_t} < \zeta t \big) + 
\textrm{P}( \nu_m > t \, , \, T_{\tilde{N}_t} \ge \zeta t ) \, .
\end{equation}

%\bigskip \bigskip

We first deal with the second term in the right hand side of (\ref{ineqct}). Before $\nu_m$ the Brownian motion cannot hit $0$. 
Thus $\textrm{P}( \nu_m > t \, , \, T_{\tilde{N}_t} \ge \zeta t )$ is bounded by the probability that the Brownian motion starting 
at $m$ hits 0 after time $\zeta t$ which is an $O(m/\sqrt{t})$. 

%\bigskip \bigskip

From now on, we only consider the first term on the right hand side of (\ref{ineqct}). 

First note that for each $j=0,m$, $(N^j_t)$ stochastically dominates a renewal process whose renewal times are distributed as the 
square root of an exponential law with parameter
\begin{equation}\label{an}
 a_n := c \, \Big( 1 + \frac{d}{n} \Big) \, (1+ \tau)^{-4} \, ,
\end{equation}
truncated at $d_n^\prime \log n$. Indeed, for $t\leq d_n^\prime \log n$, we have that
$$
P(S^j_i > t) \leq 
\textrm{P} \left(\left \{\bigcup_{u \in [0,t]} \left[ y \pm c \, \frac{1 + d/n}{1 + \tau} \, u \right] \times \{Z^j_{i-1}+u\}\right\} 
\bigcap \mathcal{P}^\star_n = \emptyset \right) \, . 
$$
The right hand side is bounded above by the probability that 
$$
\bigcup_{u \in [0,t]} \left[ y \pm c \, \frac{1 + d/n}{1 + \tau} \, u \right] \times \{u\}
$$
contains no point of a Poisson point process of parameter $(1+\tau)^{-3}$. This is readily checked to equal 
$$
e^{- \, a_n t^2} \, .
$$

From the Large Deviations Principle for the Law of Large Numbers of renewal processes, see Theorem 2.3 and Lemma 2.6 in \cite{tiefeng}, 
we have that, for $j=0,m$, $P(N^j_t \leq bt)$ decays exponentially
fast for any fixed constant $b$ smaller than 
$$
\Big( \int_0^{\infty} e^{- \, 2c\, (1+ \tau)^{-4} \,u^2} du \Big)^{-1} \, ,
$$  
which is the inverse of the mean waiting time of the renewal process. 
Writing
$T_{\tilde{N}_t} = T_{N^0_t+ N^m_t}$, the previous result implies that
$$
\textrm{P}\big(\nu_m > t, T_{\tilde{N}_t} < \zeta t \big)
$$
is bounded above by a term decaying exponentially in $n$ plus 
\begin{equation}
\label{eq:ineqct2}
\textrm{P}\big(\nu_m > t,  T_{\left\lfloor bt \right\rfloor} < \zeta t\big) \, .
\end{equation}
%\begin{equation}
%\label{eq:ineqct2}
%\textrm{P}\big( T_{\left\lfloor bt \right\rfloor} < \min\{\zeta t,\tilde{\nu}\} \big) \, .
%\end{equation}
%where $\tilde{\nu}$ is the first meeting time of two independent Brownian motions. 
From now on, we follow the arguments presented in \cite{cfd,v}. 
However, we need to make sure that the time increments $T_i - T_{i-1}$ 
are sufficiently large outside an event with probability $O(1/\sqrt{t})$. 

\medskip

Let $W_0 = T_0=0$ and $W_i = T_{i} - T_{i-1}$, $i\geq1$. We now write
$$
T_{\left\lfloor bt \right\rfloor}= \displaystyle\sum_{i=1}^{{\left\lfloor bt \right\rfloor}} ( T_{i} - T_{i-1} ) = \displaystyle\sum_{i=1}^{\left\lfloor bt \right\rfloor} W_{i} \, .
$$
Fix $\gamma > 0$ and $\tilde{b} < b$, then the probability in (\ref{eq:ineqct2}) is bounded above by 
\begin{eqnarray}
\label{eq:ineqct3}
%\lefteqn{
&P\big( \# \{ 1\leq i \leq \left\lfloor bt \right\rfloor : 0 < |Z^m_{i-1}| \leq \gamma \} \ge \lfloor (b - \tilde{b}) t \rfloor \big)&%}
\nonumber \\
& + P\big(  \# \{ 1\leq i \leq \left\lfloor bt \right\rfloor : |Z^m_{i-1}| > \gamma \} \ge \lfloor \tilde{b} t \rfloor, 
T_{\left\lfloor bt \right\rfloor} <\zeta t \big). &
\end{eqnarray}

\medskip

\begin{lem} 
\label{coalclose}
There exists $\gamma > 0$ sufficiently small such that for every $\tilde{b} < b$, there exist $\alpha = \alpha(\tilde{b},b) > 0$ and $\beta = \beta(\tilde{b},b) <\infty$ such that
$$
P\big( \# \{ 1\leq i \leq \left\lfloor bt \right\rfloor : 0 < |Z^m_{i-1}| \leq \gamma \} \ge \lfloor (b - \tilde{b}) t \rfloor \big) 
\leq \beta e^{-\alpha t} \, .
$$
\end{lem}

\medskip

The proof of Lemma \ref{coalclose} is postponed to the end of this section.

\medskip

Let $\cn=\{ 1\leq i <\infty : |Z^m_{i-1}| > \gamma \} $,
and $F_t = \big\{ \# [\cn\cap(0, bt )]\ge \lfloor \tilde{b} t \rfloor \big\}$. 
The second term in (\ref{eq:ineqct3}) is then bounded above by
\begin{equation}
\label{eq:ineqct4}
P\left(  F_t, {\sum_{i \leq \left\lfloor bt \right\rfloor : |Z^m_{i-1}| > \gamma}} W_i < \zeta t \right) 
\leq 
e^{\zeta t} \,\textrm{E} \left[ 1_{{F}_t} \,
\prod_{i \leq \left\lfloor bt \right\rfloor :\, |Z^m_{i-1}| > \gamma}e^{- W_i}\right] \, . 
\end{equation}
%$$
%E_t = \left\{ \, \max\big\{ S^j_i : j=0,m, \ 1\leq j \leq \left\lfloor bt \right\rfloor \big\} %\, < \, \frac{3}{2} \log t \, \right\} 
%$$
%Using standard estimates for the distribution of the maximum of iid random variables, the definition of the system and properties of the Poisson point process, we verify that $P(F_t^c)$ is of order $O(1/\sqrt{t})$ uniformly in $n$. Therefore, it remains to estimate
%\begin{equation}
%\label{eq:ineqct4}
%P \Big( F_t \, , \, \sum_{j \leq \left\lfloor bt \right\rfloor : |Z^m_i| > \gamma} W_j < \min\{\zeta t,\tilde{\tau}\} 
%\Big) \, .
%\end{equation}

Let $I_1<I_2<\ldots$ denote the elements of $\cn$ in increasing order. Then the expectation on the right of~(\ref{eq:ineqct4})
is bounded above by
\begin{eqnarray}
\nonumber
&\ds{\sum_{i_1,\ldots,i_{\lf \tilde bt\rf}\leq bt}
\textrm{E} \Big[ \prod_{k=1}^{\lfloor \tilde{b} t \rfloor }  e^{- W_{i_k}},\, 
I_1=i_1, ... , I_{\lfloor \tilde{b} t \rfloor } =i_{\lfloor \tilde{b} t \rfloor } \Big] }&\\
%\nonumber
&{\ds\leq\ce_{\lf \tilde bt\rf}\sum_{i_1,\ldots,i_{\lf \tilde bt\rf-1}\leq bt}
\textrm{E} \Big[ \prod_{k=1}^{\lfloor \tilde{b} t \rfloor -1}  e^{- W_{i_k}},\, 
I_1=i_1, ... , I_{\lfloor \tilde{b} t \rfloor -1} =i_{\lfloor \tilde{b} t \rfloor -1} \Big] 
\label{eq:zeta1}
\leq\cdots\leq \prod_{k=1}^{\lfloor \tilde{b} t \rfloor }\ce_k,}&
\end{eqnarray}
where, for $k\leq1$, $\ce_k=\sup\textrm{E} \Big[ e^{-  W_{i_k}} \Big| \mathcal{G}_{T_{i_k-1}} \Big]$,
with the sup taken over histories up to $T_{i_k-1}$ such that $|Z^m_{i_k-1}|>\gamma$.

Therefore to finish the proof we need to show that $\ce_k$
%\begin{equation}
%\label{eq:expeW}
% \textrm{E} \Big[ e^{-  W_{i_k}} \Big| \mathcal{G}_{T_{i_k-1}} \, , \, i_1 < ... < i_{k} \in {F}_t \Big]
%\end{equation}
is uniformly bounded way from one. To do this, we need some information on the distribution of the random vectors $(U_{i_k},V_{i_k})$.

\medskip

\begin{lem} \label{p00}
For every $0 < \delta < \frac{\gamma}{8c(1+\tau)}$ sufficiently small, there exists $\vartheta = \vartheta(\gamma,\delta) < 1$ such that
\begin{equation}
\label{eq:c_1}
\sup_{M>\gamma} P\big( \, \min\{ U_i,V_i \} \leq \delta \, \big| \, \mathcal{G}_{T_{i-1}} \, , |Z^m_{i-1}| = M \big) < \vartheta.
\end{equation}
\end{lem}

\medskip

The proof of Lemma \ref{p00} is postponed to the end of this section.

\medskip

Let $W_{-\delta,\delta}$ be the exit time of interval $(-\delta,\delta)$ by a standard Brownian motion. 
By Lemma \ref{p00}, we have that $\ce_k$
%the expectation in (\ref{eq:expeW}) 
is bounded above by the supremum over $M > \gamma$ of
\begin{eqnarray}
\label{eq:zeta2}
\lefteqn{\!\!\!\!\!\!\!\!\!\!\!\!\!\!\!\!\!\!\!\!\!\!\!\!\!\!\!\!\!\!\!\!\!\!\!\!\!\!\!\!\!\!
\textrm{E}\left[ e^{-W_{i_k}} \big| \min\{ U_{i_k},V_{i_k} \} > \delta \, , |Z^m_{i_{k}-1}| = M \right] \textrm{P}\left( \min\{ U_{i_k},V_{i_k} \} > \delta \big| |Z^m_{i_{k}-1}| = M  \right) + }  \nonumber \\
& & + \ \textrm{P}\left( \min\{ U_{i_k},V_{i_k} \} \leq \delta \big| |Z^m_{i_{k}-1}| = M  \right)  \nonumber \\
& \leq & \textrm{E}\left[ e^{-W_{-\delta,\delta}} \right] \left(1-\textrm{P}\left( \min\{ U_{i_k},V_{i_k} \} \leq \delta \big| |Z^m_{i_{k}-1}| = M  \right) \right) \nonumber \\
& & + \ \textrm{P}\left( \min\{ U_{i_k},V_{i_k} \} \leq \delta \big| |Z^m_{i_{k}-1}| = M  \right)  \nonumber \\
& = & \textrm{P}\left( \min\{ U_{i_k},V_{i_k} \} \leq \delta \big| |Z^m_{i_{k}-1}| = M  \right) (1- \tilde{\vartheta}) + \tilde{\vartheta} \nonumber \\
& \leq & \vartheta (1- \tilde{\vartheta}) + \tilde{\vartheta}
\end{eqnarray}
where $\vartheta < 1$ is given by Lemma \ref{p00} and $\tilde{\vartheta} = \textrm{E}\left[ e^{- W_{-\delta,\delta}} \right] < 1$. 
Now, choose $\zeta$ such that $\beta = e^{\zeta} [\vartheta (1- \tilde{\vartheta}) + \tilde{\vartheta}] < 1$, then from (\ref{eq:zeta1}) and (\ref{eq:zeta2}) we have that (\ref{eq:ineqct4}) is bounded by $\beta^t$. This completes the proof. $\square$

%\bigskip 

\bigskip

The rest of this section is devoted to the proofs of Lemmas \ref{coalclose} and \ref{p00}, used in the proof of 
Proposition~\ref{prop:coaltime1}.

%\bigskip 

\medskip

\noindent \textbf{Proof of Lemma \ref{coalclose}} 

The proof is based in the intuitively clear fact that if $\bar X^m_{\tilde{S}_i}$ is close to 
$\bar X^0_{\tilde{S}_{i}}$ then $\bar X^m_{s}$ and $\bar X^0_{s}$ 
should coalesce with high probability at time $\tilde{S}_{i+1}$.
%
%In order to make an argument out of this, we need some control on the height of the cones defining the processes $bar X^m_{s}$ and $bar X^0_{s}$. 
Let 
$$
J = \{ 1\leq i \leq  bt  : 0 < |Z^m_i| \leq \gamma \} \, .
$$
We will estimate
\begin{equation}
\label{eq:JJ1J2}
P\big( \# J  \ge \lf\hat b t\rf \big) \, ,
\end{equation}
where $\hat b=b-\tilde b$.

Let $I'_1<I'_2<\ldots$ denote the elements of $\{ 1\leq i <\infty  : 0 < |Z^m_i| \leq \gamma\}$ in incresing order.
For $i\geq1$, $j=0,m$, let $Y_i=j$ if the $i$-th transition point of $(\tilde N_t)_{t\geq0}$ is a transition point of $(N^j_t)_{t\geq0}$.
(See paragraph of~(\ref{Njt}).) Then the latter probability is bounded above by

\begin{eqnarray}
\nonumber
&\ds{\sum_{{i_1,\ldots,i_{\lf \hat bt\rf}\leq bt}\atop{\ell_1,\ldots,\ell_{\lf \hat bt\rf}=0\mbox{\tiny { or }}m}}
\textrm{P} \Big[  \bigcap_{k=1}^{\lfloor \hat{b} t \rfloor }  \{Z_{i_{k+1}} \neq 0, I'_k=i_k, Y_k=\ell_k\}\Big] }&\\
%\nonumber
&{\ds\leq\ce'_{\lf \hat bt\rf}
\sum_{{i_1,\ldots,i_{\lf \hat bt\rf-1}\leq bt}\atop{\ell_1,\ldots,\ell_{\lf \hat bt\rf-1}=0\mbox{\tiny { or }}m}}
\textrm{P} \Big[  \bigcap_{k=1}^{\lfloor \hat{b} t \rfloor-1}  \{Z_{i_{k+1}} \neq 0, I'_k=i_k, Y_k=\ell_k\}\Big]
\label{eq:l9}
\leq\cdots\leq \prod_{k=1}^{\lfloor \hat{b} t \rfloor }\ce'_k,}&
\end{eqnarray}
where, for $k\leq1$, $\ce'_k=\sup\textrm{P} \big[ Z_{i_{k+1}} \neq 0 \big| \mathcal{F}_{\tilde S_{i_k}} \big]$,
with the sup taken over $x,y\in\R$, $0<y-x\leq\gamma$, $\tilde s\geq0$, $\ell=0,m$, and histories up to 
$\tilde S_{i_k}=\tilde s$ 
such that $\bar X^0_{\tilde s}=x$, $\bar X^m_{\tilde s}=y$, $Y_{i_k}=\ell$.
Here ${\cal F}_t$ is the $\s$-algebra generated by $\cp_n^\star$ restricted to $\R\times[0,t]$, $t\geq0$.

\begin{figure}
\begin{center}
\input{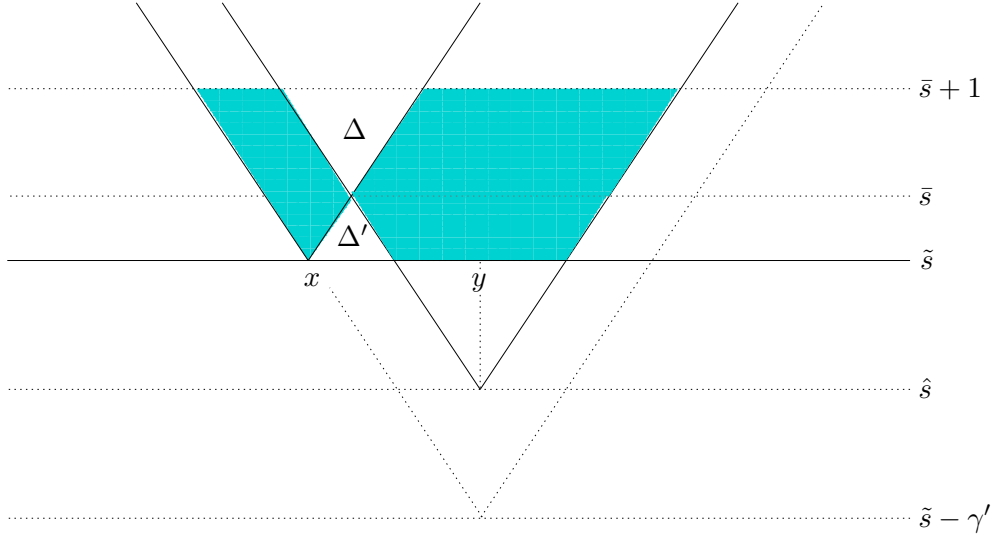}
\end{center}
\caption{$\ct$ and $\hat\ct$ are the inverted triangles with lower vertices at $(x,\tilde s)$ and $(y,\hat s)$, respectively, 
(with top portions missing); $\gamma'$ is such that the left hand boundary of $\ct''_{(y,\tilde s-\gamma')}$ touches 
$(x,\tilde s)$.}
\label{fig:coalregion}
\end{figure}

It is enough now to argue that $\ce'_k$ is bounded below away from 1 uniformly in $k$ and all large $n$. 
With $x$, $y$, $\tilde s$
fixed as above and $\ell=0$ (the case $\ell=m$ is essentially the same), let $\ct=\ct''_{(x,\tilde s)}$
and $\hat\ct=\ct''_{(y,\hat s)}$, where $\hat s$ is the last transition point of $(N^m_t)_{t\geq0}$
before $\tilde s$ (see~(\ref{eq:ttt}) and~(\ref{t_alt}) for the definition of $\ct''_{(y,s)}$).
See Figure \ref{fig:coalregion}. 

Given that $\bar X^m_{\tilde s}>\bar X^0_{\tilde s}$ and that $Y_{i_k}=0$, we must have that 
$\hat\ct\cap\{\R\times[0,\tilde s]\}$ does not contain $(\bar X^0_{\tilde s},s)$
--- otherwise, since the interior of $\hat\ct\cap\{\R\times[0,\tilde s]\}$ must contain no point of $\cp_n^\star$,
we would have that $\bar X^m_{\tilde s}=\bar X^0_{\tilde s}$.
All this implies that $\tilde s-\gamma'<\hat s<\tilde s$, where $\gamma'=\gamma'(y-x)$ is the 
sup over $z>0$ such that $\ct''_{(y,\tilde s-z)}\cap\{\R\times[0,\tilde s]\}$ 
does not contain $(\bar X^0_{\tilde s},\tilde s)$. 
See Figure \ref{fig:coalregion}. 

Let now $\bar s$ denote the time coordinate of the point where the right hand boundary of $\ct$ meets the 
left hand boundary of $\hat\ct$ (see Figure \ref{fig:coalregion}). We define the following regions of $\H$.
Let $Q$ be the quadrangle bounded by the horizontal lines at times $\tilde s$, $\bar s+1$,
the left hand boundary of $\ct$ and the right hand boundary of $\hat\ct$. Let $\Delta$ be the triangle 
bounded by the horizontal line at time $\bar s+1$, the right hand boundary of $\ct$ and the left hand boundary of 
$\hat\ct$. And let $\Delta'$ be the triangle
bounded by the horizontal line at time $\tilde s$, the right hand boundary of $\ct$ and the left hand boundary of 
$\hat\ct$. Let $Q'=Q\setminus\{\Delta\cup\Delta'\}$. See Figure \ref{fig:coalregion}, where $Q'$ appears
shaded. 

Let $A_k=A_k(x,y,\tilde s,\hat s)$ denote the event that $Q'\cap\cp_n^\star=\emptyset$ and 
$\Delta\cap\cp_n^\star\neq\emptyset$.
It is a simple yet tedious matter to verify that the area of $\Delta$ is bounded away from zero, and the 
area of $Q'$ is bounded as $x,y,\tilde s,\hat s$ vary within their restricted range (notice that 
$\gamma'$ is bounded above by $\gamma'(\gamma)<\infty$ within that range). This readily implies 
that $P(A_k)$ is bounded away from zero within that range.

Since $\ce'_k\leq\inf P(A_k^c)$, with the inf taken over $x,y,\tilde s,\hat s$ varying within their restricted range, the 
result follows from the uniform positive bound on $P(A_k)$ and~(\ref{eq:l9}).
 $\square$

\bigskip 

%\medskip

\noindent \textbf{Proof of Lemma \ref{p00}} 

%Let us start with a simple observation. By standard properties of Poisson Point Processes, 
%if both processes $\bar X^0_s$ and $\bar X^m_s$ jump at the same time, then they coalesce. 
%If they coalesce by time $\tilde{S}_{i}$, then 
%$|Z^m_{i}|= 0$ and $U_{i} = - Z^m_{i-1} < - \gamma < - \delta$. With this in mind, 

\begin{figure}
\begin{center}
\input{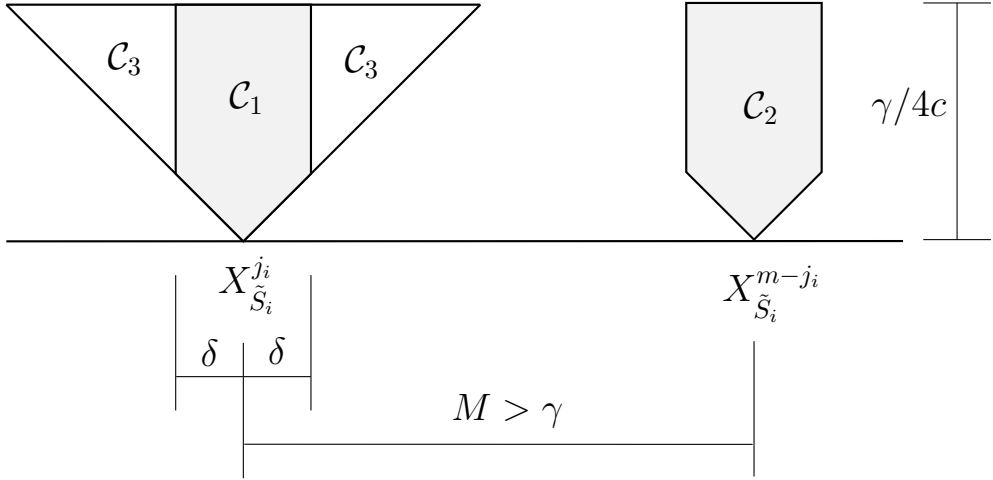}
\end{center}
 \caption{Representation of $\mathcal{C}_1$, $\mathcal{C}_2$ and $\mathcal{C}_3$}
\label{fig:coalUV}
\end{figure}

As in the proof of Lemma~\ref{coalclose}, let us again suppose that $Y_i=0$ (and again, the case $Y_i=m$ is
argued similarly).
Given 
$\mathcal{G}_{T_{i-1}}$ and $|Z^m_{i-1}| = M$, $M > \gamma$, the probability of the event 
$\{\min\{ U_{i},V_{i} \} > \delta\}$ is bounded below by the probability of the event 
$F_1 \cap F_2 \cap F_3$, where
$$
F_1 = \{ \mathcal{C}_1 \cap \mathcal{P}^{\star}_n = \emptyset \} \, , \ F_2 = \{ \mathcal{C}_2 \cap \mathcal{P}^{\star}_n = \emptyset \} 
\, \textrm{ and } \,
F_3 = \{ \mathcal{C}_3 \cap \mathcal{P}^{\star}_n \neq \emptyset \} \, ,
$$
with
$$
\mathcal{C}_1 = \mathcal{T}_{(\bar X^{0}_{\tilde{S}_{i-1}},\tilde{S}_{i-1}),\frac{\gamma}{4c}}^{\prime \prime} 
\cap 
\left\{\left[ \bar X^{0}_{\tilde{S}_{i-1}} \pm \delta \right] 
\times \left[\tilde{S}_{i-1}, \tilde{S}_{i-1} + \frac{\gamma}{4c}\right]\right\} \, ,
$$
$$
\mathcal{C}_2 = \mathcal{T}_{(\bar X^{m}_{\tilde{S}_{i-1}},\tilde{S}_{i-1}),\frac{\gamma}{4c}}^{\prime \prime} 
\cap 
\left\{\left[ \bar X^{m}_{\tilde{S}_{i-1}} \pm \delta \right] 
\times \left[\tilde{S}_{i-1}, \tilde{S}_{i-1} + \frac{\gamma}{4c}\right]\right\} \, ,
$$
and
$$
\mathcal{C}_3 = \mathcal{T}_{(\bar X^{0}_{\tilde{S}_{i-1}},\tilde{S}_{i-1}),\frac{\gamma}{4c}}^{\prime \prime}\setminus\mathcal{C}_1 \, .
$$
See Figure \ref{fig:coalUV}.
Since $\mathcal{C}_1$,  $\mathcal{C}_2$ and  $\mathcal{C}_3$ are disjoint with 
Area$(\mathcal{C}_1) =\mbox{Area}(\mathcal{C}_2) \leq \frac{\delta \, \gamma}{2c}$ and 
$$
\mbox{Area}(\mathcal{C}_3) = 
\mbox{Area} \big( \mathcal{T}_{(\bar X^{0}_{\tilde{S}_{i-1}},\tilde{S}_{i-1}),\frac{\gamma}{4c}}^{\prime \prime} \big) 
- \mbox{Area}(\mathcal{C}_1) 
\ge \frac{\gamma^2}{16c^2(1+\tau)} - \frac{\delta \gamma}{2 c}
= \frac{\gamma}{2c} \big( \frac{\gamma}{8 c (1+\tau)} - \delta \big) > 0 \, .
$$
Thus,
\begin{eqnarray*}
&P\big( \, \min\{ U_i,V_i \} \ge \delta \, \big| \, \mathcal{G}_{T_{i-1}} \, , |Z^m_{i-1}| = M \big) 
\ge  P\big( \, F_1 \cap F_2 \cap F_3 \, \big| \, \mathcal{G}_{T_{i-1}} \, , |Z^m_{i-1}| = M \big) &\\
& =  P\big( F_1 \big) \, P\big( F_2 \big) \, P\big( F_3 \big)
\, \ge \, e^{-\frac{\delta \gamma}{c}} \big( 1 - e^{-\frac{\gamma}{2c} \big( \frac{\gamma}{8 c (1+\tau)} - \delta \big)} \big) \, .&
\end{eqnarray*}
$\square$

\medskip

\begin{rmk}
 \label{rmk:ext}
 We will make use below of two extensions of Proposition~\ref{prop:coaltime1}, whose proofs are straightforward 
 adaptations of or require minimal additional arguments with respect to the above proof. One extension is for 
 an independent pair $((\hat X^0_t)_{t\geq0},(\hat X^m_{t})_{t\geq0})$, with marginals equally distributed to 
 those of $((X^0_t)_{t\geq0},(X^m_{t})_{t\geq0})$. The other is as follows.
 \begin{equation}
\label{eq:coaltime-ext}
\textrm{P} (\nu_m > t|\mathcal{F}_s) \leq \frac{C \, |X^m_{s}-X^0_{s}|}{\sqrt{t-s}} \, , \, \textrm{ for every } t>s>0.
\end{equation}
\end{rmk}

%%%%%%%%%%%%%%%%%%%%%%%%%%%%%%%%%%%%%%%%%%%%%%%%%%%%%%%%%%%%%%%
\section{Verification of condition I} 
%%%%%%%%%%%%%%%%%%%%%%%%%%%%%%%%%%%%%%%%%%%%%%%%%%%%%%%%%%%%%%%
\label{sec:I}

\subsection{Convergence of a single trajectory of $ \tilde\G'''_n$}\label{ssingle}

In this subsection we consider a single trajectory starting from a deterministic point of
$\bar\La'_n$ of the form $({\sqrt n}z_0,nt_0)$ for given fixed $(z_0,t_0)\in(-\infty,\infty)\times[0,\tau)$.
We will denote it as above by $\hat\g'''_{\mbox{}_{{\sqrt n}z_0,nt_0}}$. And we will denote its
rescaled version $D(\hat\g'''_{\mbox{}_{{\sqrt n}z_0,nt_0}})$ by $\tilde\g'''_{\mbox{}_{z_0,t_0}}$.

We will establish a weak convergence result of $\tilde\g'''_{\mbox{}_{z_0,t_0}}$ to a Brownian motion
as follows.

\begin{prop}\label{single}
As $n\to\infty$
 \begin{equation}
  \label{eq:single}
  \tilde\g'''_{\mbox{}_{z_0,t_0}}\Rightarrow B_{\mbox{}_{z_0,t_0}},
 \end{equation}
 a Brownian trajectory with diffusion coefficient 
  \begin{equation}
  \label{eq:om}
  \o=\frac1{\sqrt{6}\,\pi^{1/4}}\,(\tan\t)^{3/4} \, ,
 \end{equation}
 starting at $(z_0,t_0)$, and ending at time $\tau$,
 where $\Rightarrow$ in this case means convergence in distribution in the uniform topology on continuous
 trajectories from $[0,\tau]\to\R$.
\end{prop}

%\begin{rmk}\label{rmk:path}
% Notice that the portion of $\hat\g'''_{\mbox{}_{{\sqrt n}z_0,nt_0}}$ from the second vertex on is a 
%path of $\hat\G''_n$.
%\end{rmk}

\begin{rmk}\label{rmk:edge}
 As reasoned in the proof of Lemma~\ref{lm:vdist}, it is not very important how we join the (directed) edges of 
 $\tilde\g'''_{\mbox{}_{z_0,t_0}}$, provided they stay within the isosceles triangle described towards the end
 of that proof, which, as argued in that proof, is the case of the edges of the trajectories of $\hat\G''_n$, and
 is clearly also the case for linear interpolations of successive vertices. So, it is also the case for the 
 edges of $\tilde\g'''_{\mbox{}_{z_0,t_0}}$. Below we will consider a jump version of 
 $\tilde\g'''_{\mbox{}_{z_0,t_0}}$, for which the same also holds.
\end{rmk}

%\begin{rmk}\label{rmk:path}
% Notice that with high probability the vertices of $\hat\g'''_{\mbox{}_{{\sqrt n}z_0,nt_0}}$ from the second 
%vertex on coincide with the vertices of the path of $\hat\G''_n$ starting at this second vertex. From the
%arguments of Remark~\ref{rmk:edge} and the proof of Lemma~\ref{lm:vdist} the Hausdorff distance between 
%the reescaled versions of those paths (in the appropriate space) vanishes as $n\to\infty$.
%\end{rmk}

\noindent{\bf Proof} 

By the horizontal translation invariance of the model we may take $z_0=0$. We will for simplicity also take
$t_0=0$. The argument for other cases is an easy adaptation.

Let us consider the {\em jump} version of $\hat\g'''_0$ defined as follows. 
\begin{equation}\label{eq:zpt}
Z'_t=\sum_{i=1}^J(Y'_i-Y'_{i-1})\,1_{\{t\geq S'_i\}},
\end{equation}
where $(Y'_i,S'_i)=\hat s_i'''(0)$. By Lemma~\ref{lm:idd}, $(Y'_i,S'_i)_{i\geq1}$ is distributed
like $(Y_i,S_i)_{i\geq1}$. So it is enough to show the convergence to $B_0$ of $\{Z^{(n)}_t,\,0\leq t\leq\tau \}$,
where $Z^{(n)}_t=\frac1{\sqrt n}Z_{tn}$, and 
\begin{equation}\label{eq:zt}
Z_t=\sum_{i\geq1}X_i\,1_{\{t\geq S_i\}},
\end{equation}
where $X_i=Y_i-Y_{i-1}$.
%with ${\cal J}=\inf\{i\geq1:\,S_i>\tau n\}.$

%Let 
%\begin{equation}\label{eq:nt}
%N_t=\sum_{i\geq1}1\{t\geq S_i\}.
%\end{equation}

%Then, $Z_t=Y_{N_t}$. 

We start by establishing a law of large numbers for $S_{rn}$.

\paragraph{Law of large numbers for $S_{rn}$.}

\begin{lem}\label{lm:lln}
Given $J>0$, we have that almost surely
\begin{equation}
\label{eq:lln}
\sup_{0\leq r\leq \frac Jn}\le|\frac1 nS_{\lf rn\rf}-\frac {\hat cr}{1-\hat cr}\ri|\to0
\end{equation} 
as $n\to\infty$, where $\hat c$ is a positive constant to be defined below.
\end{lem}

\noindent{\bf Proof of Lemma~\ref{lm:lln}} 

It is convenient to go back to $S'_{rn}$ instead, and use the map back to $\La_n$, where the
issue involves essentialy iid increments, rather than location dependent ones.

Indeed, let us recall that $S'_i=\hat s_i'''(0)_2=\frac{n-|\hat s_i(0,-n)|}{|\hat s_i(0,-n)|/n}$ (where $\hat s_i(x)$
starting from a deterministic point of $\La_n$ is defined as in the beginning of Subsection~\ref{pre},
using $\ct_{x,\log n}$). We readily find that 
\begin{equation}
\label{eq:l1}
n-|\hat s_i(0,-n)|=:S''_i=\sum_{j=1}^i R_i,
\end{equation} 
where given $S''_{i-1}$, $R_i$ is distributed as
\begin{equation}
\label{eq:l2}
\P(R_i>u|S''_{i-1})=\exp\{-\mbox{area of $\ct_{(0,S''_{i-1}),u}$}\}1_{\{u<\log n\}}.
\end{equation} 

One may readily check from our discussions on Subsection~\ref{pre} (see e.g.~Figure~\ref{fig:tt}
and~(\ref{eq:a1})) that for $i=1,\ldots, \hat I$ the area of $\ct_{(0,S''_{i-1}),u}$ is bounded 
from below and from above by respectively $cu^2$ and $(c+c'/n)u^2$, where $c'$ is a constant.
(Recall that with high probability $\hat I$ is the first $i$ for which $|\hat s_i(0,-n)|<\a n$ and
$\hat I=J$.)
We conclude that we may dominate $R_1,R_2,\ldots$ from above and from below by iid sequences of random 
variables $R'_1,\ldots,R'_{\hat I}$ and $R''_1,\ldots,R''_{\hat I}$, respectively, where
\begin{eqnarray}
\label{eq:l3}
\P(R'_1>u)&=&e^{-cu^2}1_{\{0<u<\log n\}},\\
\label{eq:l4}
\P(R''_1>u)&=&e^{-(c+c'/n)u^2}1_{\{0<u<\log n\}}.
\end{eqnarray} 

It follows readily that we may with probability larger than $1-e^{-c''(\log n)^2}$ replace
$R'_1,\ldots,R'_{\hat I}$ and $R''_1,\ldots,R''_{\hat I}$ by respectively  
$\hat R'_1,\ldots,\hat R'_{\hat I}$ and $\hat R''_1,\ldots,\hat R''_{\hat I}$, independent 
random variables such that
\begin{eqnarray}
\label{eq:l5}
\P(\hat R'_i>u)&=&e^{-cu^2},\\
\label{eq:l6}
\P(\hat R''_i>u)&=&e^{-(c+c'/n)u^2},
\end{eqnarray} 
where $c''$ is a positive constant.

It now follows from standard large deviation estimates that outside an event of exponentially 
small probability in $n$, and inside the event of the previous paragraph, given $\a'\in(0,\a)$,
$J$ is smaller than $\frac1{\hat c}(1-\a')n$. Again applying standard large 
deviation estimates we get that 
\begin{equation}
\label{eq:llna}
\sup_{0\leq r\leq \frac{\hat I}n}\le|\frac1 nS''_{\lf rn\rf}-\hat cr\ri|\to0
\end{equation} 
as $n\to\infty$, using the fact that $\E(\hat R'_1)=\hat c:=\frac12\sqrt{\pi/c}$ and 
$\E(\hat R''_1)=\frac12\sqrt{\pi/(c+c'/n)}$.

The result follows from the representation of $S'$ in terms of $S''$ discussed at the beginning of 
this proof.
$\square$

\medskip

We continue with the proof of Proposition~\ref{single}.

\paragraph{Convergence of finite dimensional distributions.}

Let $k\geq1$ and $0<t_1<\ldots<t_k\leq\tau$. Given $S_1,S_2\ldots$ satisfying~(\ref{eq:lln}), which is
an event of full measure, we have that the increments of $Z^{(n)}$ are independent. Let us consider 
$Z^{(n)}_{t_2}-Z^{(n)}_{t_1}$. We write it as
\begin{equation}
 \label{eq:f1}
\frac1{\sqrt n}\sum_{i\geq1}X_i\,1_{\le\{t_1<\frac{S_i}n\leq t_2\ri\}}.
\end{equation}
Taking the log of the Laplace transform of the above random variable, conditional on $S_1,S_2\ldots$, we get
\begin{equation}
 \label{eq:f2}
 \sum_{i\geq1}\log\k\le(\frac{\l c_n}{\sqrt n}\frac{T_i1_{\{T_i<L_n\}}}{1+\frac{S_{i-1}}n}
 1_{\le\{t_1<\frac{S_i}n\leq t_2\ri\}}\ri),
\end{equation}
where $\k(x)=\sinh(x)/x$, $\l$ is the argument of the transform, and $T_1,T_2,\ldots$ are independent,
with $T_i$ distributed as~(\ref{eq:idf}), with $s=S_{i-1}$.

Since $\k(x)=1+\frac16x^2+O(x^4)$, we may estimate~(\ref{eq:f2}) by
\begin{equation}
 \label{eq:f3}
\frac{\l^2 c_n^2}{6n} \sum_{i\geq1}\frac{T^2_i1_{\{T_i<L_n\}}}{\le(1+\frac{S_{i-1}}n\ri)^2}
 1_{\le\{t_1<\frac{S_i}n\leq t_2\ri\}}+%\ri)+
 \mbox{const }\frac{(\log n)^4}{n^2}\sum_{i\geq1}1_{\le\{t_1<\frac{S_i}n\leq t_2\ri\}}.
\end{equation}

The estimates in the proof of Lemma~\ref{lm:lln} imply that the second term of the above sum %right of~(\ref{eq:f3}) 
is almost surely negligible as $n\to\infty$. Let us analyse the first term.

Given $0<\eps<t_1$, Lemma~\ref{lm:lln} can be applied to get that the first term in (\ref{eq:f3}) is
almost surely bounded from above and below respectively by
\begin{equation}
 \label{eq:f4}
(1\pm\eps)^2\frac{\l^2 c_n^2}{6n} 
\sum_{i=\frac n{\hat c}\frac{t_1\pm\eps}{1+t_1}}^{\frac n{\hat c}\frac{t_2\pm\eps}{1+t_2}}
\frac{\tilde T^2_i}{\le(1-\hat c\frac in\ri)^2} 
%\frac{T^2_i}{\le(1+\frac{S_{i-1}}n\ri)^4}
\end{equation}
for all large $n$, where $\tilde T_i=T_i/\le(1+S_{i-1}/n\ri)^2$. 
(That the indicator $1_{\{T_i<L_n\}}$ may be dropped follows from the fact that 
$\ct''_{\hat s''_i(z_0,t_0),L_n}\cap\cp'\ne\emptyset$ for all $i=0,1,\ldots,J$,
with high probability, which can be argued as in the proof of Lemma~\ref{lm:id} above.)

\begin{rmk}\label{rmk:tpi}
 In order to estimate latter expression, let us first observe that from~(\ref{eq:idf})
 %and Lemma~\ref{lm:lln}, 
 we may dominate the distribution of $\tilde T_i$, $i=1,2,\ldots$, 
 above and below by $\check T_i$, $i=1,2,\ldots$, iid random variables such that
 \begin{equation}\label{eq:idfp}
\P(\check T_1>v)=e^{-\check c_nv^2}%,\quad\P(T''_1>v)=e^{-c_n''v^2},\,v>0.
\end{equation}
and $\check c_n\to \hat{c}$ as $n\to\infty$, 
where $(\check c_n)$ does not depend on $S$ and varies as upper and lower bounds).
\end{rmk}

Let us analyse thus
\begin{equation}
 \label{eq:f5}
(1\pm\eps)^2\frac{\l^2 c_n^2}{6n} 
\sum_{i=\frac n{\hat c}\frac{t_1\pm\eps}{1+t_1}}^{\frac n{\hat c}\frac{t_2\pm\eps}{1+t_2}}
\frac{\check T^2_i}{\le(1-\hat c\frac in\ri)^2}.
%\frac{T^2_i}{\le(1+\frac{S_{i-1}}n\ri)^4}
\end{equation}

A standard large deviation estimate tells us that the latter expression is bounded above and below respectively
by
\begin{equation}
 \label{eq:f7}
(1\pm\eps)^2\frac{\l^2 c_n^2}{6\check c_n} 
\sum_{i=\frac n{\hat c}\frac{t_1\pm\eps}{1+t_1}}^{\frac n{\hat c}\frac{t_2\pm\eps}{1+t_2}}
%\sum_{i=2cn\frac{t_1\pm\eps}{1+t_1\pm\eps}}^{2cn\frac{t_2\pm\eps}{1+t_2\pm\eps}}
\frac{1/n}{\le(1-\hat c\frac in\ri)^2}\pm\eps,
%\frac{\E(T^2_i|S_{i-1})}{\le(1+\frac{S_{i-1}}n\ri)^4}.
\end{equation}
also for all $n$ sufficiently large, where we have used the fact, as follows from~(\ref{eq:idfp}), that $\E(\check T^2_1)=1/\check c_n$.

Now, the latter sum is a Riemann sum for the integral
\begin{equation}
 \label{eq:f8}
\int_{\frac1{\hat c}\frac{t_1\pm\eps}{1+t_1}}^{\frac1{\hat c}\frac{t_2\pm\eps}{1+t_2}}
\frac{dx}{\le(1-\hat cx\ri)^2}=
\frac1{\hat c}\int_{\frac{t_1\pm\eps}{1+t_1}}^{\frac{t_2\pm\eps}{1+t_2}}
\frac{dy}{\le(1-y\ri)^2}.
\end{equation}

Since $\eps$ is arbitrary, we find that~(\ref{eq:f5}) converges almost surely as $n\to\infty$
to
\begin{equation}
 \label{eq:f9}
\frac{\l^2 c}{6\hat c}\int_{\frac{t_1}{1+t_1}}^{\frac{t_2}{1+t_2}}
\frac{dy}{\le(1-y\ri)^2}=\l^2\frac{c}{6\hat c}(t_2-t_1)=:\l^2\o^2(t_2-t_1) \, ,
\end{equation}
where
$$
  \o=\sqrt{\frac{c}{6\hat c}}=\frac1{\sqrt{6}\,\pi^{1/4}}\,(\tan\t)^{3/4}.
$$
Collecting the above steps, we readily conclude that
given $S_1,S_2\ldots$ in a set of full measure, the increments of $Z^{(n)}$ 
converge to independent Gaussian random variables with variance given by $\o^2$ times
the time increments. This establishes the convergence of the finite dimensional distributions
of $Z^{(n)}$ given $S_1,S_2\ldots$ in a set of full measure to those of Brownian motion with 
diffusion coefficient $\o$. 
%The result follows by integrating on $S_1,S_2\ldots$. 

Let us now check tightness of $Z^{(n)}$ given $S_1,S_2\ldots$  in a set of full measure. 
Along with the finite dimensional distribution convergence result, that implies that
the statement of Proposition~\ref{single} holds for the distribution of
$Z^{(n)}$ given $S_1,S_2\ldots$ in a set of full measure.
The (unconditional) result then follows by integration. $\square$

\paragraph{Tightness.} 

We will verify standard tightness criteria for the distribution of
$Z^{(n)}$ given $S_1,S_2\ldots$ in a set of full measure. 
We may assume $S_1,S_2\ldots$ satisfies~(\ref{eq:lln}).

Given the convergence of finite dimensional distributions established above, it is enough to
verify condition  {\bf (b)} of Corollary 7.4, page 129 in Ethier and Kurtz (1986). For that
it is enough to show that given $\eps,\d>0$ and $\cj=\lc \d^{-1}\tau\rc$, making 
$t_j=j\d$, $j=0,1,\ldots,\cj$, we have that 
$$
\lim_{\d\to0}\limsup_{n\to\infty}\sum_{j=1}^\cj
P\le(\sup_{s,t\in[t_{j-1},t_j]}|Z^{(n)}_t-Z^{(n)}_s|>\eps\ri)=0,
$$
where $P(\cdot)=\P(\cdot|S_1,S_2\ldots)$. Indeed, we are going to show that
\begin{equation}
 \label{eq:t1}
\lim_{\d\to0} \delta^{-1} \limsup_{n\to\infty} \sup_{1 \leq j \leq \cj}
P\le(\sup_{s,t\in[t_{j-1},t_j]}|Z^{(n)}_t-Z^{(n)}_s|>\eps\ri)=0 \, .
\end{equation}

It is enough to get this result replacing $\sup_{s,t\in[t_{j-1},t_j]}|Z^{(n)}_t-Z^{(n)}_s|$
by $\sup_{t\in[t_{j-1},t_j]}|Z^{(n)}_t-Z^{(n)}_{t_{j-1}}|$.

Since 
\begin{equation}
  \label{eq:t2}
Z^{(n)}_t-Z^{(n)}_{t_{i-1}}=\frac1{\sqrt n}\sum_{i\geq1}X_i\,1_{\le\{t_{j-1}<\frac{S_i}n\leq t\ri\}}
=:\frac1{\sqrt n}W_t,
\end{equation}
and using Markov's inequality, we get that
\begin{equation}
 \label{eq:t3}
%\lim_{\d\to0}\limsup_{n\to\infty}\sum_{j=1}^\cj
P\le(\sup_{t\in[t_{j-1},t_j]}|Z^{(n)}_t-Z^{(n)}_{t_{j-1}}|>\eps\ri)=
P\le(\sup_{t\in[t_{j-1},t_j]}|W_t|>\eps{\sqrt n}\ri)\leq\frac1{\eps^4n^2}E(M_j^4)\leq\frac{\mbox{const}}{n^2}E(W_{t_j}^4),
\end{equation}
where $M_j=\sup_{t\in[t_{j-1},t_j]}|W_t|$, and we have used the $L^p$ {\em maximum inequality}, valid here since
under $P$, the $Y_i$'s are independent and have zero mean.

Now, $E(W_{t_j}^4)/n^2$ is equal to
\begin{equation}
 \label{eq:t4}
\frac1{n^2}\sum_{i\geq1}(E(X_i^4)-E^2(X_i^2))\,1_{\le\{t_{j-1}<\frac{S_i}n\leq t_j\ri\}}+
            \le(\frac1n\sum_{i\geq1}E(X_i^2)\,1_{\le\{t_{j-1}<\frac{S_i}n\leq t_j\ri\}}\ri)^2 \, ,
\end{equation}

The first term of (\ref{eq:t4}) is positive. Dropping
$E^2(X_i^2)$ and using~(\ref{eq:xys}), we find that it is bounded above by constant times
\begin{equation}
 \label{eq:t5}
\frac1{n^2}\sum_{i\geq1}T_i^4\,1_{\le\{t_{j-1}<\frac{S_i}n\leq t_j\ri\}},
\end{equation}
and arguing as in the estimation of~(\ref{eq:f4}), we find that the latter sum is of order $n$ outside 
an event of $\P$-probability exponentially small. Thus the first term of~(\ref{eq:t4})
is almost surely negligible as $n\to\infty$. 
The squared term on~(\ref{eq:t4}) may likewise be upper bounded by
\begin{equation}
 \label{eq:t6}
\frac1n\sum_{i\geq1}T_i^2\,1_{\le\{t_{j-1}<\frac{S_i}n\leq t_j\ri\}},
\end{equation}
and again an argument like the one to estimate~(\ref{eq:f4}) yields an almost sure upper bound for the $n$ limit
of~(\ref{eq:t6}) of constant times $\d$. 

Substituting successively in~(\ref{eq:t4}) and~(\ref{eq:t3}), we get~(\ref{eq:t1}).
$\square$

\medskip

\begin{rmk}
 \label{rmk:ext_ip}
 Below we will need extensions of Proposition~\ref{single} to the case of conditional distributions of
 $\tilde\g'''_{\mbox{}_{z_0,t_0}}$ and its jump version given the history up to a deterministic or stopping
 time. These follow by virtually the same reasoning as above, with minimal, straifgtforward modifications.
\end{rmk}

\bigskip 

\subsection{Convergence of a finite number of trajectories of $\tilde\G'''_n$}

Now that we have convergence of single trajectories to Brownian Motion, we can prove condition $I$ following the same steps of the 
proof presented in \cite{cv}. A similar approach is undertaken also in~\cite{cfd}, so we will allow ourselves to be somewhat sketchy
in our arguments for this subsection.
%It is a consequence of the next proposition. As in Section \ref{sec:coaltime}, we consider the 
%extension of the system of coalescing trajectories evolving in $\mathbb{R}^2$. 
With respect to~\cite{cv}, here we have the advantage that trajectories cannot 
cross each other and the disadvantage that the trajectories are not evolving according to a discrete space-time lattice.

\medskip

\begin{prop}
\label{prop:coalBM}
Let $(z_0,t_0)$, $(z_1,t_1)$, ... , $(z_m,t_m)$ be $m+1$ distinct points in $\mathbb{R}\times[0,\tau)$. Then
$$
\Big( \tilde{\gamma}^{\prime \prime \prime}_{z_0,t_0} , ... , \tilde{\gamma}^{\prime \prime \prime}_{z_m,t_m} \Big) 
\Longrightarrow^D \big( B_{z_0,t_0} , ... , B_{z_m,t_m} \big) \, ,
$$
where $B_{z_0,t_0}$, ... , $B_{z_m,t_m}$ are coalescing Brownian Motions with constant diffusion coefficient $\o$ 
starting at $(z_0,t_0)$,...,$(z_m,t_m)$. 
\end{prop}

We prove Proposition \ref{prop:coalBM} by induction on $m$, the case $m=0$ having been treated in Proposition~\ref{single}.
We may start by supposing that $t_j < t_m$, $j=1,...,m-1$. 
From the induction hypothesis, conditioning on the history up to $t_m$, we may indeed reduce to the case where
$t_0=t_1=\cdots=t_m$, and, relabeling if necesssary, $z_0<z_1<\cdots<z_m$. 

%\medskip

%We already have the convergence of $ \tilde{\gamma}^{\prime \prime \prime}_{z_k,t_k}$ to a 
%Brownian Motion starting at $z_k$ at time $t_k$. 

%\medskip

Let us fix a uniformly continuous bounded function $H:D([t_m,\tau])^{m+1} \rightarrow \mathbb{R}$, with the following property.
Let us start by defining {\em coalescence} operators of two trajectories % $C:D([t_m,\tau])^{2} \rightarrow D([t_m,\tau])^{2}$,
as follows.

Given $\g,\g'\in D([t_m,\tau])$ such that $\g(0)<\g'(0)$, let $\check t=\check t(\g,\g')=\inf\{t\in[t_m,\tau]:\,\g(t)\geq\g'(t)\}$, 
with $\inf\emptyset=\infty$. 
Now let $C(\g,\g')=(\g,\check\g)$, with $\check\g(t)=\g'(t)$ for $t<\check t$, and $\check\g(t)=\g(t)$ for $t\geq \check t$.
This should be seen as coalescence with the path (initially) to the left.
(The cases where $\g(0)=\g'(0)$ are immaterial for our purposes, and can be defined arbitrarily, say in such a way that
either $C(\g,\g')=(\g,\g)$ or $C(\g,\g')=(\g',\g')$.)

Let us now define a coalescence operator of $m+1$ trajectories. 
Suppose $\g_0,\ldots,\g_m\in D([t_m,\tau])$ such that $\g_0(t_m)<\cdots<\g_m(t_m)$. Then 
$C_m(\g_0,\ldots,\g_m)=(\check\g_0,\ldots,\check\g_m)$, where $\check\g_{0}=\g_{0}$, and for $k=1,\ldots,m$, 
$(\check\g_{k-1},\check\g_{k})=C(\check\g_{k-1},\g_{k})$. %, and $(\check\g_{m-1},\check\g_{m})=C^l(\g_{m-1},\g_{m})$. 
We may say that under $C_m$, $\g_{0}$ remains invariant, and the paths $\g_1,\ldots,\g_{m}$ coalesce to the left.

Now for the above mentioned property of $H$. We require $H$ to be invariant under coalescence, 
in the following sense. Given $\g_0,\ldots,\g_m$ as above, we ask that
$H(\g_0,\ldots,\g_m)=H\circ C_m(\g_0,\ldots,\g_m)$.

%Supposing that $(\tilde{\gamma}^{\prime \prime \prime}_{z_0,t_0} , ... , \tilde{\gamma}^{\prime \prime \prime}_{z_{m-1},t_{m-1}})$ 
%converges in distribution to 
%$(B_{z_0,t_0} , ... , B_{z_{m-1},t_{m-1}})$, w

We will show that
\begin{eqnarray}
\label{I2}
\lim_{n \rightarrow \infty} 
|\textrm{E}[H(\bar{\gamma}^{\prime \prime \prime}_{z_0,t_0} , ... , \bar{\gamma}^{\prime \prime \prime}_{z_m,t_m})] 
- \textrm{E}[H(B_{z_0,t_0} , ... , B_{z_m,t_m})]| = 0 \, ,
\end{eqnarray}
where $\bar{\gamma}^{\prime \prime \prime}_{z_k,t_k}$ is the jump version of $\tilde{\gamma}^{\prime \prime \prime}_{z_k,t_k}$
(similarly as above). Notice that $C_m$ is almost surely continuous with respect to the product Wiener measure on $D([t_m,\tau])^{m+1}$
(under the sup norm).
By induction and the definition of convergence in distribution, we obtain Proposition \ref{prop:coalBM}.

\medskip

We start by taking a version of $\bar{\gamma}^{\prime \prime \prime}_{z_m,t_m}$ which is independent of 
$(\bar{\gamma}^{\prime \prime \prime}_{z_0,t_0} , ... , \bar{\gamma}^{\prime \prime \prime}_{z_{m-1},t_{m-1}})$.
Let
$$
\nu = \tau \wedge \inf \big\{ s \ge 0 : \, ( \bar{\gamma}^{\prime \prime \prime}_{z_m,t_m} (s) 
     - \bar{\gamma}^{\prime \prime \prime}_{z_{m-1},t_{m-1}}(s) ) \leq n^{-\frac{1}{8}} \big\}.
$$
For every $n$, let $\mathcal{P}^{\star \star}_n$ be a Poisson point process which is also independent of 
$\mathcal{P}^{\star}_n$ and has the same intensity measure given in (\ref{intensityPstar}). 
Let 
$\mathcal{Q}^{\star}_n=\{\mathcal{P}^{\star}_n\cap\{\R\times[0,\nu]\}\}\cup\{\mathcal{P}^{\star \star}_n\cap\{\R\times(\nu,\tau]\}\}$.
One readily checks that $\mathcal{Q}^{\star}_n$ is equally distributed with $\mathcal{P}^{\star}_n$. Now let 
$\bar{\gamma}^{\star}_{z_m,t_m}$
be the as the path $\bar{\gamma}^{\prime \prime \prime}_{z_m,t_m}$, except that using $\mathcal{Q}^{\star}_n$
rather than $\mathcal{P}^{\star}_n$. It may be checked that for all large enough $n$, 
$\bar{\gamma}^{\star}_{z_m,t_m}$ is independent of 
$(\bar{\gamma}^{\prime \prime \prime}_{z_0,t_0} , ... , \bar{\gamma}^{\prime \prime \prime}_{z_{m-1},t_{m-1}})$.
Notice that $\bar{\gamma}^{\star}_{z_m,t_m}$ equals $\bar{\gamma}^{\prime \prime \prime}_{z_m,t_m}$ up to time $\nu$.

\medskip
We are now ready to prove (\ref{I2}). The expression inside the lim sign there is bounded above by
\begin{eqnarray}
\label{I1}
& & \qquad |\textrm{E}[H(\bar{\gamma}^{\prime \prime \prime}_{z_0,t_0} , ... , \bar{\gamma}^{\prime \prime \prime}_{z_m,t_m})] 
- \textrm{E}[H(\bar{\gamma}^{\prime \prime \prime}_{z_0,t_0} , ... ,
                \bar{\gamma}^{\prime \prime \prime}_{z_{m-1},t_{m-1}}, \bar{\gamma}^{\star}_{z_m,t_m})]| \nn \\
& & \qquad \quad + | \textrm{E}[H(\bar{\gamma}^{\prime \prime \prime}_{z_0,t_0} , ... ,
                \bar{\gamma}^{\prime \prime \prime}_{z_{m-1},t_{m-1}}, \bar{\gamma}^{\star}_{z_m,t_m})]
- \textrm{E}[H(B_{z_0,t_0} , ... , B_{z_m,t_m})]|. 
\end{eqnarray}

By the induction hypothesis and Proposition~\ref{prop:coaltime1}, we have that the second term in (\ref{I1}) goes to zero 
as $n$ goes to $+\infty$. 
So we only have to deal with the first term in (\ref{I1}). This is bounded above by
\begin{eqnarray}
\label{I4}
\textrm{E}\big[\big|H(\bar{\gamma}^{\prime \prime \prime}_{z_0,t_0} , ... , \bar{\gamma}^{\prime \prime \prime}_{z_m,t_m}) 
- H(\bar{\gamma}^{\prime \prime \prime}_{z_0,t_0} , ... ,
                \bar{\gamma}^{\prime \prime \prime}_{z_{m-1},t_{m-1}}, \check{\gamma}^{\star}_{z_m,t_m})\big|
\mathbb{I}_{\nu<\tau} \big],
\end{eqnarray}
where $\check{\gamma}^{\star}_{z_m,t_m}$ such that 
$(\bar{\gamma}^{\prime \prime \prime}_{z_{m-1},t_{m-1}}, \check{\gamma}^{\star}_{z_m,t_m})
=C(\bar{\gamma}^{\prime \prime \prime}_{z_{m-1},t_{m-1}}, \bar{\gamma}^{\star}_{z_m,t_m})$.

\medskip

To deal with the expectation in~(\ref{I4}), we define the coalescence times
$$
\s = \inf \{ s \ge 0 : \bar{\gamma}^{\prime \prime \prime}_{z_{m-1},t_{m-1}} = \bar{\gamma}^{\prime \prime \prime}_{z_m,t_m} \} \quad \textrm{and} \quad 
\s^\star = \inf \{ s \ge 0 : \bar{\gamma}^{\prime \prime \prime}_{z_{m-1},t_{m-1}}  \geq \bar{\gamma}^{\star}_{z_m,t_m} \} \, .
$$
The times $\tau$ and $\tau^\star$ have the tail of their distributions $O(1/\sqrt{tn})$ --- see Proposition~\ref{prop:coaltime1} and
Remark~\ref{rmk:ext}. 
%

%Note that on $\mathcal{A}^c_{n,\tau}$ we have $\nu < \tau$. 
Define the event
\begin{eqnarray*}
\mathcal{C}_{n,\tau}\= \Big\{ \sup_{0 \leq s \leq \tau} |\check{\gamma}^{\star}_{z_m,t_m} (s) 
     - \bar{\gamma}^{\prime \prime \prime}_{z_m,t_m} (s)| \ge n^{-\frac{1}{16}}\log n\Big\}\\
\=  \Big\{ \sup_{\nu \leq s \leq \tau} |\check{\gamma}^{\star}_{z_m,t_m} (s) 
     - \bar{\gamma}^{\prime \prime \prime}_{z_m,t_m} (s)| \ge n^{-\frac{1}{16}} \log n \Big\},
\end{eqnarray*}
where the second equality follows from the fact that $\bar{\gamma}^{\star}_{z_m,t_m}$ equals 
$\bar{\gamma}^{\prime \prime \prime}_{z_m,t_m}$ up to time $\nu$.

Now
$
\textrm{P} \Big( \mathcal{C}_{n,\tau}, \nu<\tau \Big) 
$
is bounded above by
\begin{equation}
\label{I3}
\textrm{P} \Big( \mathcal{C}_{n,\tau},\,\nu<\tau,\,\{\s, \, \s^\star \in [\nu , \nu + n^{-\frac{1}{8}}]\} \Big) 
+ \textrm{P} \Big(\s > \nu + n^{-\frac{1}{8}}\Big) + \textrm{P} \Big(\s^\star > \nu + n^{-\frac{1}{8}} \Big) .
\end{equation}
By Proposition \ref{prop:coaltime1} and its extensions --- see Remark~\ref{rmk:ext} ---, the latter two terms in (\ref{I3}) 
are bounded above by $2 \,\frac{n^{\frac{3}{8}}}{n^{\frac{7}{16}}} = 2 \, n^{-\frac{1}{16}}$. 
On the other hand, since 
$\check{\gamma}^{\star}_{z_m,t_m}=\bar{\gamma}^{\prime \prime \prime}_{z_m,t_m}$ %=\bar{\gamma}^{\prime \prime \prime}_{z_{m-1},t_{m-1}}$
after max$\{\s,\s^\star\}$, the first term in (\ref{I3}) is bounded above by
$$
\textrm{P} \Big(\sup_{\nu \leq s \leq (\nu + n^{-\frac{1}{8}}) \wedge \tau } 
|\check{\gamma}^{\star}_{z_m,t_m} (s) - \bar{\gamma}^{\prime \prime \prime}_{z_m,t_m} (s)| 
\ge n^{-\frac{1}{16}}\log n    ,\, \nu<\tau \Big), 
$$
and this is in turn bounded above by
\begin{eqnarray*}
 &\ds{\textrm{P} \Big( \sup_{0 \leq s \leq n^{\frac{7}{8}}} 
\frac{\big|\breve{\gamma}^{\prime \prime \prime}_{z_m,t_m} (\nu+s)-\breve{\gamma}^{\prime \prime \prime}_{z_m,t_m} (\nu)\big|}
{n^{\frac{7}{16}}} \ge \frac{\log n}{2} \,\Big|\,\nu\Big)}&\\
 &\ds{+\textrm{P} \Big( \sup_{0 \leq s \leq n^{\frac{7}{8}}} 
\frac{\big|\breve{\gamma}^{\star}_{z_m,t_m} (\nu+s)-\breve{\gamma}^{\star}_{z_m,t_m} (\nu)\big|}
{n^{\frac{7}{16}}} \ge \frac{\log n}{2} \,\Big|\,\nu\Big)},&
\end{eqnarray*}
where $\breve{\gamma}^{\prime \prime \prime}_{z_m,t_m}$ is the unscaled version of $\bar{\gamma}^{\prime \prime \prime}_{z_m,t_m}$
(also the jump version of $\hat{\gamma}^{\prime \prime \prime}_{z_m,t_m}$), and likewise for
$\breve{\gamma}^{\star}_{z_m,t_m}$ with respect to $\bar{\gamma}^{\star}_{z_m,t_m}$.
By Proposition~\ref{single} and its extension (see Remark~\ref{rmk:ext_ip}), the first probability above 
goes to zero as $n$ goes to infinity, and by a virtually forthright 
extension of those results for $\breve{\gamma}^{\star}_{z_m,t_m}$, so does the second probability.
 
Finally we have that (\ref{I4}) is bounded above by a term that converges to zero as $n \rightarrow \infty$ plus
\begin{eqnarray}
 \textrm{E}\big[\big|H(\bar{\gamma}^{\prime \prime \prime}_{z_0,t_0} , ... , \bar{\gamma}^{\prime \prime \prime}_{z_m,t_m})] 
- H(\bar{\gamma}^{\prime \prime \prime}_{z_0,t_0} , ... ,
                \bar{\gamma}^{\prime \prime \prime}_{z_{m-1},t_{m-1}}, \check{\gamma}^{\star}_{z_m,t_m})\big|
\mathbb{I}_{\mathcal{C}^c_{n,\tau},\nu<\tau} \big].
\end{eqnarray}
By the uniform continuity of $H$ the rightmost expectation in the previous expression converges to zero as $n$ goes to $+\infty$.
\bigskip

%%%%%%%%%%%%%%%%%%%%%%%%%%%%%%%%%%%%%%%%%%%%%%%%%%%%%%%%%%%%%%%
\section{Verification of conditions $B_1$ and $E$} 
%%%%%%%%%%%%%%%%%%%%%%%%%%%%%%%%%%%%%%%%%%%%%%%%%%%%%%%%%%%%%%%
\label{sec:B1E}

%\bigskip \bigskip

\subsection{Verification of condition $B_1$}

\bigskip

By spatial translation invariance of the system we have to show that
$$
\limsup_{\epsilon \rightarrow 0+} \limsup_{n \rightarrow +\infty} 
\sup_{ t > \beta } \sup_{ t_0 \in \mathbb{R} } \textrm{P} ( \eta_{\tilde{\G}'''_n}(t_0,t;0,\epsilon) \ge 2) = 0 \, .
$$
We note that $\eta_{\tilde{\G}'''_n}(t_0,t;0,\epsilon) \geq 2$ if and only if leftmost and rightmost 
trajectories of $\tilde{\Gamma}'''_n$ crossing the interval $[0,\epsilon]$ at time $t_0$ have not met up to time $t_0+t$. 
By Proposition \ref{prop:coalBM}, this pair of trajectories converge to those of two coalescing Brownian motions starting 
at points $(0,t_0)$ and $(\epsilon,t_0)$. Then, it is straightforward to get that
$$
\limsup_{n\rightarrow
\infty}\mathbb{P}(\eta_{{n}}(t_0,t;0,\epsilon)\geq2)=
2\Phi(\epsilon/\sqrt{2t}\,) - 1 \leq 2\Phi(\epsilon/\sqrt{2 \beta}\,) - 1,
$$
where $\Phi(\cdot)$ is the standard normal distribution function. From the previous inequality we obtain $B_1$ by taking 
the limit as $\epsilon$ goes to $0$.

\bigskip \bigskip

\subsection{Verification of condition $E$}

\bigskip

We will work here with the following sets of paths
\begin{eqnarray*}
 \breve{\Gamma}'''_n=\{\breve\gamma'''_x,\,x\in\cp'\},
\end{eqnarray*}
which is the jump version of $\{\hat\gamma'''_x,\,x\in\cp'\}$. %, and its scaled version
%\begin{eqnarray*}
%\bar{\Gamma}'''_n=\{D(\gamma):\,\gamma\in\breve{\Gamma}'''_n\}.
%\end{eqnarray*}

%defined in section \ref{sec:coaltime}.  

To simplify notation we drop the triple primes in the remainder of this subsection, 
writing $\breve{\Gamma}_n$ in place of $\breve{\Gamma}'''_n$,
and $\breve{\g}_x$ in place of $\breve{\g}'''_x$. 
As before, define the set of diffusively rescaled paths of $\breve{\Gamma}_n$ as 
\begin{equation}
 \bar\G_n=\{D(\g);\,\g \in \breve\G_n\}.
\end{equation}

Notice that in $\breve{\Gamma}_n$ we only have paths starting from the points of the Poisson point process $\cp'$.
Nevertheless, it follows from arguments above that
\begin{equation}
 d_{\ch_{0}^{\tau,\tau}}(\tilde\G'''_n,\bar\G_n)\to0
\end{equation} 
with high probability (even though $\bar\G_n$  is in principle not in $\ch_{0}^{\tau,\tau}$ --- but could be included,
as c\`adl\`ag trajectories ---; see arguments in the proof of 
Lemma~\ref{lm:vdista} and Remark~\ref{rmk:edge}), and thus subsequential limits of $\tilde\G'''_n$ and $\bar\G_n$ 
coincide (along the same subsequences).
%
%
%Clearly for any $\mathcal{X'''}$ and $\mathcal{X}$ subsequential limits of respectively $\tilde{\Gamma}'''_n$ 
%and $\tilde{\Gamma}_n$, we have
%$$
%\textrm{E}[ \hat{\eta}_\mathcal{X'''} (t_0,t;a,b) ] \leq  \textrm{E}[ \hat{\eta}_\mathcal{X} (t_0,t;a,b) ] \, .
%$$ 
Then, it is enough to show that  $\bar{\Gamma}_n$ satisfies condition $E$.

\medskip

We will  follow \cite{nrs} and \cite{s} closely, with similar notation, which we now introduce. We fix $\mathcal{X}$ 
as a subsequencial limit of $\bar{\Gamma}_n$, which is a tight sequence, as follows from Proposition B.2 of~\cite{finr1}, 
since the paths of each of its elements are noncrossing, and, as seen above, converge to Brownian motions 
(see Proposition~\ref{single} and its proof above). 
For any system of space time paths $\mathcal{Y}$, given $T\in\R$,
we write set $\mathcal{Y}^{T^-}$ as the set of paths in $\mathcal{Y}$ that start at some time $s<T$. 
We also write $\mathcal{Y}(T)$ to represent the set of intersection points of all paths in $\mathcal{Y}$
with $\R\times\{T\}$. 
Note that, $\hat{\eta}_\mathcal{Y} (t_0,t;a,b) = \#(\mathcal{Y}^{t_0^-}(t_0+t) \cap (a,b))$.

\medskip

In the proof of condition $E$, the first result we need to show is that $\mathcal{X}^{t_0^-}(t_0 + \epsilon)$ 
is a locally finite point process. Here the proof is more complicated than the lattice random walk case presented 
for instance in \cite{s}. The first step is to prove that $\breve{\Gamma}^{T^-}_n(T)$ is  in a certain sense 
locally finite uniformly in $T>0$ (note that here we are considering $\breve{\Gamma}_n$ as a set of paths starting 
at space-time points in $\mathcal{P}'$). This is the object of the next result.

%\bigskip

\begin{lem}\label{gammalf} There exists a constant $C > 0$, which does not depend on the scaling parameter $n$, such that
$$ 
\textrm{E} \Big[ \# \big( \breve{\Gamma}^{T^-}_n(T) \cap [0,M) \big) \Big] \leq C \, M \, .
$$
\end{lem}

%\medskip

\noindent \textbf{Proof} 

We say that a point $(x,s) \in \{[0,M) \times [0,T)\} \cap \mathcal{P}'$ \emph{touches} $[0,M)\times\{T\}$ if the path 
$\breve{\gamma}_{x,s}$ does not meet any other point of $\mathcal{P}'$ during the time interval $[s,T]$. By the definition 
of the random paths in $\breve{\Gamma}_n$, if $\breve{\gamma}_{x,s}$ touches $[0,M)\times\{T\}$, then it is constantly 
equal to $x$ in the time interval $[s,T)$. Note that $\breve{\Gamma}^{T^-}_n(T) \cap \{[0,M)\times\{T\}\}$ is equal to the random 
set of points that touch $[0,M)\times\{T\}$. 

Now fix $L = \frac c{1+\tau} \wedge 1$. Enlarging $M$ if necessary, we can suppose that $M/L$ is an integer. 
For $j=1,...,M/L$ and $1\leq k\leq\lc T\rc$, let
$D_{j,k} = [(j-1)L,jL) \times [T-k,T-k+1)$, and let $A_{j,k}$ be the random sets of points in 
$D_{j,k} \cap \mathcal{P}'$ that touch $[0,M)\times\{T\}$, and also let
$$B_{j,k} = D_{j,k} \cap \mathcal{P}'.$$ 
For $k>\lc T\rc$, let $B_{j,k} \equiv \emptyset$.

\begin{figure}
\begin{center}
\input{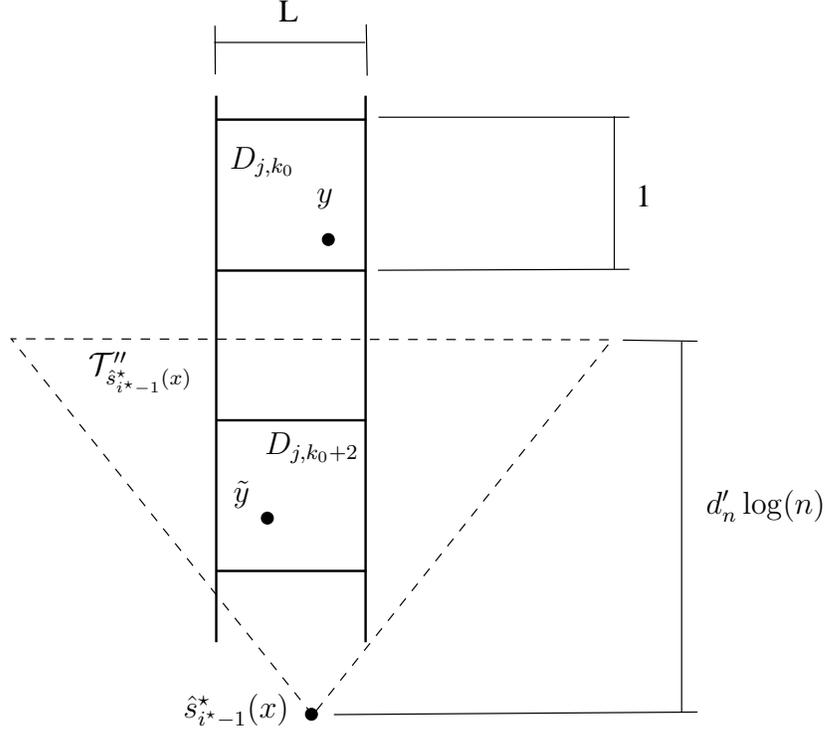}
\end{center}
 \caption{
If $y \in B_{j,k_0}$ and  $\tilde{y} \in B_{j,k_0+ 2}$ and 
$x \in  A_{j,k}$ for some $k > k_0 + 2$, then either 
$\tilde y\in\mathcal{T}_{\hat{s}''_{i^\star-1} (x), d_n^\prime \log n }^{\prime \prime}$, 
or $y\in\mathcal{T}_{\hat{s}''_{i^\star} (x), d_n^\prime \log n }^{\prime \prime}$.}
\label{fig:touch0M}
\end{figure}

We claim that $B_{j,k_0} \neq \emptyset$ and $B_{j,k_0+2} \neq \emptyset$ implies that $\# A_{j,k} = \emptyset$ for every 
$k > k_0 + 2$. To prove the claim, let $y \in B_{j,k_0}$ and  $\tilde{y} \in B_{j,k_0+ 2}$. Suppose that there exists 
$x=(x_1,x_2) \in  A_{j,k}$ for some $k > k_0 + 2$. Let $\d_0=x_2$ and, for $i\geq0$, $\d_{i+1}=\d_{i}+d_n^\prime(\d_{i}) \log n$
--- here we are making explicit that $d_n^\prime(\cdot)$ is a function; see its definition on the paragraph of~(\ref{eq:ttt}) above.
Let $i^\star$ be the largest $i\geq0$ such that $\d_{i^\star}<T-k_0+1$. Then either $\hat{s}''_{i^\star} (x)$, which has to
equal $(x_1,\d_{i^\star})$, lies in $D_{j,k_0}\cup D_{j,k_0+1}$, in which case 
$\tilde y\in\mathcal{T}_{\hat{s}''_{i^\star-1} (x), d_n^\prime \log n }^{\prime \prime}$, 
in contradiction to the fact that $x$ touches $M$, or it lies below $D_{j,k_0+1}$, in which case
$y\in\mathcal{T}_{\hat{s}''_{i^\star} (x), d_n^\prime \log n }^{\prime \prime}$, again
in contradiction to the fact that $x$ touches $M$, see Figure \ref{fig:touch0M}. And the claim is established.

Define $\beta_j = \min \big\{ k \ge 1 : B_{j,3k} \neq \emptyset \textrm{ and } B_{j,3k+2} \neq \emptyset \big\}$,
with $\min\emptyset=(\lc T\rc - 2)/3$. 
The random variables $\beta_j$, $j=1,...,M/L$, are iid random variables stochastically dominated by a geometric 
distribution of parameter $(1-e^{-L})^2$. Moreover, 
$$
\# \big( \breve{\Gamma}^{T^-}_n(T) \cap \{[0,M)\times\{T\}\} \big) 
= \sum_{j=1}^{M/L} \sum_{k=1}^{3\beta_j + 2} \# A_{j.k} 
\leq \sum_{j=1}^{M/L} \sum_{k=1}^{3 \beta_j + 2} \# B_{j.k}  \, .
$$
Now the $B$'s and $\beta$'s are not independent, but $\# B_{j.k} | \beta_j = m$ is stochastically dominated by 
$1+\zeta$ where $\zeta$ is a Poisson distribution of parameter $L$, for every $1 \leq k \leq 3m +2$. Therefore
\begin{eqnarray*}
\textrm{E} \Big[ \# \big( \breve{\Gamma}^{T^-}_n(T) \cap \{[0,M)\times\{T\}\} \big) \Big] & \leq & 
\sum_{j=1}^{M/L} \textrm{E} \Big[ \sum_{k=1}^{3 \beta_j + 2} \# B_{j.k}  \Big] \ \leq \
\frac{M}{L} \, \textrm{E} \Big[ \textrm{E} \Big[ \sum_{k=1}^{3 \beta_j + 2} \# B_{j.k} \Big| \beta_j \Big] \Big] \\
& \leq & \frac{M}{L} \, \textrm{E} \big[ 1 + \zeta \big] \, \textrm{E} \big[ {3 \beta_j + 2} \big] 
\ = \ \frac{M}{L} \, (1+L) \, \Big( \frac{3}{(1-e^{-L})^2}  + 2\Big) \, .
\end{eqnarray*}   
$\square$

\bigskip

\begin{lem}\label{density}
There exists a positive constant $C>0$, which does not depend on the scaling parameter $n$, such that 
$$ 
\textrm{E} \Big[ \# \big( \breve{\Gamma}^{T^-}_n(T+t) \cap \{[0,M)\times\{T+t\}\} \big) \Big] \leq \frac{C \, M}{\sqrt{t}} \, ,
$$
for every $M>0$.
\end{lem}

%\begin{Lem}\label{density}
%Denote by  $O(t_0,t)$ the event that there is a path of $\hat{\G}^\star_n(t_0-)$ that visits $[-1/2,1/2]$ by time $t_0+t$. Then
%$$
%P\big(O(t_0,t)\big) \leq \frac{C}{\sqrt{t}}
%$$
%for some positive constant $C$ independent $n$.
%\end{Lem}

\medskip

\begin{rmk}
Lemma \ref{density} is a version of Lemma 2.0.7 in section 2 of \cite{s}.  The latter result holds for the difference of 
two independent continuous time random walks on $\mathbb{Z}$ which is not our case. 
\end{rmk}

\medskip

\noindent \textbf{Proof of Lemma~\ref{density}} 

First, we point out that, by the additivity of $\# \big( \breve{\Gamma}^{T^-}_n(T+t) \cap \cdot \big)$ as a set function 
in the borelians of $\R\times\{T+t\}$, % of $A$,
it is enough to consider the case $M=1$. 
%Indeed, by (spatial) translation invariance
%$$
%\textrm{E} \Big[ \# \big( \breve{\Gamma}^{t_0^-}_n(t_0+t) \cap [0,M) \big) \Big] \leq \sum_{j=1}^{\left\lceil M \right\rceil}  \textrm{E} \Big[ \# \big( \hat{\Gamma}^{t_0^-}_n(t_0+t) \cap [j,j+1) \big) \Big] =  \left\lceil M \right\rceil \, \textrm{E} \Big[ \# \big( \hat{\Gamma}^{t_0^-}_n(t_0+t) \cap [0,1) \big) \Big] \, .
%$$

Below we use Proposition \ref{prop:coaltime1} and Lemma \ref{gammalf} above, as well as an adaptation of the proof of Lemma 2.7 in~\cite{nrs} 
(or Lema 2.07 in \cite{s}) that works if we properly replace the counting variables $\xi^A_t$ and the notion of nearest neighbors sites,
as follows. 

\smallskip

For every $A \subset \mathbb{Z}$, $t>0$, and $k \in \mathbb{Z}$, let $\xi^A_t(k)$ be the number of points in 
$\breve{\Gamma}_n^{T^-} (T+t) \cap \{[k,k+1)\times\{T+t\}\}$ due to paths that also visit $[l,l+1)$ at time 
$T$ for some $l \in A$. Also define $\xi^A_t = \sum_{j \in \mathbb{Z}} \xi_t^{A}(j)$ which is the number of points 
in $\breve{\Gamma}_n^{T^-} (T+t)$ also due to paths that also visit $[l,l+1)$ at time $T$ for some $l \in A$. 
With this definition, we have that 
$e_t = \textrm{E} \Big[ \# \big( \breve{\Gamma}^{T^-}_n(T+t) \cap \{[0,1)\times\{T+t\}\} \big) \Big] = \textrm{E} \big[ \xi^\mathbb{Z}_t(0) \big]$.

\smallskip

Let $B_N = \{0,1,...,N-1\}$, by the same translation invariance argument presented in Lemma 2.7 in \cite{nrs}, we obtain that
\begin{eqnarray}
\label{eq:xiBM}
e_t \, N & = & \textrm{E} \Big[ \sum_{j=0}^{N-1} \xi^\mathbb{Z}_t(j)\Big] 
\, \leq \, \sum_{k \in \mathbb{Z}} \textrm{E} \Big[ \sum_{j=0}^{N-1} \xi^{B_N + kN}_t(j)\Big] \nonumber \\
& = &  \sum_{k \in \mathbb{Z}} \textrm{E} \Big[ \sum_{j=0}^{N-1} \xi^{B_N}_t(j+kN)\Big]  \leq \textrm{E} \Big[{\xi}_t^{B_N} \Big].
\end{eqnarray}

By Lemma \ref{gammalf}, we have that $\textrm{E} \big[ {\xi}_t^{B_N} \big]$ is finite. For every $k \in B_N$, let 
$E_k = \breve{\G}_n^{T -}(T) \cap \{[k,k+1)\times\{T\}\}$ and $\beta_k=\#E_k$. To avoid the event that $\beta_k \equiv 0$, 
we perform the enlargement $J_N = \{ 1/2 , 3/2 , ... , (2N+1)/2 \} \cup \big( \cup_{j=0}^{N-1} E_j \big)$. Thus $J_N$ is a set 
of cardinality $\beta = N + \sum_{j=0}^{N-1} \beta_k$ with nearest neighbor points at distance smaller or equal than one 
from each other. Given $\beta$, ${\xi}_t^{B_N}$ is bounded above by $\beta$ times the number of nearest neighbor pairs 
in $J_N$ that have coalesced by time $t$. Thus, by Proposition \ref{prop:coaltime1}, there exists $C_1 > 0$ such that
\begin{equation*}
\textrm{E} \big[ {\xi}_t^{B_N} \big| \beta \, \big] \leq \beta - (\beta - 1) \textrm{P}(\nu_1 \leq t) 
\leq 1 + \beta \, \frac{C_1}{\sqrt{t}} \, .
\end{equation*}
Also by Lemma \ref{gammalf}, we have that the expectation of $\beta_k$ is bounded above by $C_ 2$ for some $C_2 > 0$. 
Hence, for $C = C_1 (C_2+1)$
$$
\textrm{E} \big[ \hat{\xi}_t^{B_N} \big] \leq 1 + \frac{C \, N}{\sqrt{t}} \, .
$$
From the previous inequality , (\ref{eq:xiBM}) and the fact that $N$ is arbitrarily chosen, we have that $e_t \leq C/ \sqrt{t}$.
$\square$

\bigskip

As a straighforward application of the previous lemma, we get the following result.

\medskip

\begin{cor}\label{cor:density}
For every $t_0$, $t$, $a$, $b \in \mathbb{R}$ with $t>0$ and $a<b$, we have 
$$ 
\limsup_{n \rightarrow \infty} \textrm{E} \Big[ \hat{\eta}_{\bar{\Gamma}_n} (t_0,t;a,b) \Big] \leq \frac{C \, (b-a)}{\sqrt{t}} \, .
$$
\end{cor}

\bigskip

The previous corollary is a version of Lemma 3.5.4 in \cite{s}. From it we are able to obtain version of Lemmas 3.5.2 and 3.5.3, 
which we write  as follows.

\medskip

\begin{lem}\label{locfin} For every $t>0$,
$\mathcal{X}^{t_0^-}(t_0+t)$ is almost surely locally finite.
\end{lem}

\medskip

\begin{lem}\label{locfinconv} For every $t>0$,
$\mathcal{X}^{t_0^-}_{t_0+t}$, i.e. the set of paths starting at $\mathcal{X}^{t_0^-}(t_0+t)$ truncated before time $t_0+t$, 
is distributed as coalescing Brownian motions starting at the random set $\mathcal{X}^{t_0^-}(t_0+t)$.
\end{lem}

\bigskip

With these results, we may conclude, as in \cite{nrs} and \cite{s}, as follows. For every $\epsilon < t/2$, we have that
$$
\textrm{E}[ \hat{\eta}_\mathcal{X} (t_0,t;a,b) ] \leq
\textrm{E} \Big[ \hat{\eta}_{\mathcal{X}^{(t_0+\epsilon)^-}_{t_0+t}} (t_0+\epsilon,t;a,b) \Big] \, .
$$ 
From Lemma \ref{locfinconv}, the right hand side in the previous inequality is bounded above by 
$$
\textrm{E}[ \hat{\eta}_{\cw_0} (t_0+\epsilon,t;a,b) ] = \frac{b-a}{\sqrt{\pi (t-\epsilon)}} \, .
$$
Letting $\epsilon \rightarrow 0$, we obtain condition $E$. 

\bigskip
 
\noindent \textbf{Proof of Lemma \ref{locfin}} 

The proof is entirely analogous to that of Lemma 3.5.2 in \cite{s}. 
Indeed, all we need there and here are simple facts about weak convergence and Corollary \ref{cor:density}.

\bigskip

\noindent \textbf{Proof of Lemma \ref{locfinconv}} 

The random set $\bar\G_n^{t_0^-}(t_0+t)$ converges in 
distribution to $\mathcal{X}^{t_0^-}(t_0+t)$. Since we already have condition $I$, upon attempting to follow the proof of 
Lemma 3.5.3 in \cite{s}, we realize that all we need is a version of Lemma 3.5.5 in that paper that could 
be applied to our case, see also Remark 3.5.1. The technical drawback here is that we cannot consider $\bar\G_n^{t_0^-}(t_0+t)$ as starting 
points of random walks due to the non-Markovian property of the paths. Indeed this is the only difficulty. 
At first sight it may appear that Lemma 3.5.5 in \cite{s} holds for discrete spatial lattices only, but 
this hypothesis is not used in the proof. Indeed, it may be readily 
checked that that proof holds for random sets of $\mathbb{R}^2$ that are almost surely locally finite.
%as follows.

%To deal with the problem mentioned above, 
We continue by replacing $\bar\G_n^{t_0^-}(t_0+t)$ by a suitable  diffusively rescaled subset of 
$\mathcal{P}'$. Let $\tilde{\mathcal{P}} = \{ D(x,s) : (x,s) \in \mathcal{P}' \}$. 
For each point $x$ in a given realization of $\bar\G_n^{t_0^-}(t_0+t)$, we have that 
$\tilde{\mathcal{P}}\cap \big( \{x\} \times (0,t_0+t) \big)$ is a unitary set almost surely; 
we call its single point the ancestor of $x$, denoted by $a(x)$. Let $\bar A_n$ be the random set of ancestors of 
$\bar\G_n^{t_0^-}(t_0+t)$. We claim that $\bar A_n$ converges in distribution to $\mathcal{X}^{t_0^-}(t_0+t)$. 
Since $\bar A_n$ consists of starting points of the trajectories in $\bar{\Gamma}_n$, if we prove the claim then 
we can use Lemma 3.5.5 and adapt the proof of Lemma 3.5.3 in \cite{s} to our case. 

\smallskip

The remainder of this proof will be devoted to prove that $\bar A_n$ converges in distribution to $\mathcal{X}^{t_0^-}(t_0+t)$. 
Indeed, we show that for each fixed $M>0$ the Hausdorff distance between $\bar{A}_{M,N} := \{x=(x_1,x_2)\in \bar A_n:\,x_1\in (-M, M)\}$ 
and $\bar\G_{M,n} := \bar\G_n^{t_0^-}(t_0+t) \cap \{(-M,M)\times\{T+t\}\}$ converges to zero in probability. We denote the Hausdorff 
distance between sets in $\mathbb{R}^2$ by $\rho_H$. In order to avoid complications with notation due to scaling, we introduce the 
following the following objects. Let
$\breve\G_{M,n} := \breve\G_n^{(nt_0)^-}(n(t_0+t)) \cap (-M\sqrt n,M\sqrt n)\times\{n(t_0+t)\}$, and let $\breve A_{M,n}$ be the set 
of ancestors of $\breve\G_{M,n}$, where the definition of ancestor in is analogous to the one above but uses $\mathcal{P}'$ in place 
of $\tilde{\mathcal{P}}$. 
%The choice 
%was made to avoid complications with notation due to scaling. 
We then have that $\bar\G_{M,n}=D(\breve\G_{M,n})$, $\bar{A}_{M,n} =D(\breve A_{M,n})$, and 
\begin{equation}
\label{eq:convsets}
 \rho_H \big( \bar{A}_{M,n} , \bar\G_{M,n} \big) 
 %\leq n^{-2} \rho_H \big( \hat{A}_{M,n} , \hat\G_{M,n} \big) = 
 \leq n^{-1} \sup \Big\{ | x - a(x) | : x \in  \breve\G_{M,n} \Big\} \, .
\end{equation}

Now, from the proof of Lemma \ref{gammalf}, we have that for all $x \in  \hat\G_{M,n}$, $a(x)$ is a point that touches 
$(-M\sqrt n,M\sqrt n)$ at time $n(t_0 + t)$. Proceeding as in that proof, we make a partition of the interval $(-M\sqrt n,M\sqrt n)$ in 
$\left\lceil 2M\sqrt n/L \right\rceil$ intervals, $I_j$, of size at most $L$. We associate to each 
$I_j$, $1\leq j \leq \left\lceil 2Mn/L \right\rceil$, a random variable $\beta_j$ also as in the proof of Lemma~\ref{gammalf}. 
Then, recall that the definition of $\beta_j$ implies that no point at distance greater than $3 \beta_j + 2$. Therefore,  
for every $\epsilon > 0$, 
$$
P \big( \rho_H \big( \bar{A}_{M,n} , \bar\G_{M,n} \big) \ge \epsilon \big) \leq
P \Big( \max_{1\leq j \leq \frac{M \sqrt n}{L} }\frac{3 \beta_j + 2}{n}  \ge \epsilon  \Big).
$$

Since, as argued in the proof of Lemma \ref{gammalf} above, the random variables $\beta_j$ are iid and stochastically 
dominated by a geometric distribution, a standard argument shows that the latter probability vanishes as $n\to\infty$.
%Replacing the maximum by a sum in the last probability, using classical Markov inequality and the fact that , 
%it is straightforward to verify that the probability is of order $1/n$. 
This proves the claim. $\square$

\end{document}